\documentclass{article}
\input tesis.sty
\usepackage{amsfonts}
\usepackage{amssymb}
\usepackage{a4wide}
\usepackage[english]{babel}
\newcommand{\Pol}{{\cal P}}
\newcommand{\Har}{{\cal H}}
\newcommand{\Mon}{{\cal M}}

\newcommand{\comment}[1]{}

\newcommand{\Span}{\mbox{span}}

\renewcommand{\S}{{\mathbb S}}
\newcommand{\mS}{{\mathbb S}}
\renewcommand{\cong}{\simeq}
\newcommand{\phar}{P_H}
\newcommand{\pmon}{P_M}
\newcommand{\Bp}{B_+}
\newcommand{\Bm}{B_-}
\newcommand{\Bn}{B_0}
\newcommand{\Fm}{F_-}
\newcommand{\Fp}{F_+}
\newcommand{\gosp}{\mathfrak{osp}}
\newcommand{\Ve}{I}
\newcommand{\We}{J}

\newcommand{\adj}{\ast}
\newcommand{\End}{\mbox{End}}
\newcommand{\Aut}{\mbox{Aut}}
\newcommand{\SO}{\mbox{SO}}
\newcommand{\inpro}[2]{ \langle #1,#2 \rangle}
\newcommand{\bin}[2]
{\left(\! \!
\begin{array}{c}
#1\\
#2
\end{array}
\! \! \right)
}
\newcommand{\bina}[2]
{\left\{\! \!
\begin{array}{c}
#1\\
#2
\end{array}
\! \! \right\}
}

\begin{document}
\title{Orthogonal basis for spherical monogenics by step two branching}
\author{R.~L\'avi\v cka\thanks{e-mail:lavicka@karlin.mff.cuni.cz, 
Mathematical Institute,
Charles University,  
Sokolovsk\'a 83, 
186 75 Praha, Czech Republic}, \,\,\, V.~Sou\v cek\thanks{e-mail:soucek@karlin.mff.cuni.cz, 
Mathematical Institute,
Charles University,  
Sokolovsk\'a 83, 
186 75 Praha, Czech Republic}, \,\,\, P.~Van Lancker\thanks{e-mail:Peter.VanLancker@hogent.be, 
Department of Engineering Sciences,
University College Ghent, Member of Ghent University Association, 
Schoonmeersstraat 52, 
9000 GENT (Belgium)}}
%\date{}
\maketitle

\begin{abstract}
Spherical monogenics can be regarded as a basic tool for the study of harmonic analysis of the Dirac operator in Euclidean space $\R^m$. They play a similar role as spherical harmonics do in case of harmonic analysis of the Laplace operator on $\R^m$.
%The space of Clifford-valued spherical monogenics $ \Mon_k(\R,\C_m)$ in $\R^m$ can be regarded as a model for the irreducible representation of $\Spin(m)$ with %weight $(k+\frac{1}{2},\frac{1}{2},\ldots,\frac{1}{2})$. 
Fix the direct sum $\R^m=\R^p \oplus \R^q$. In this paper we will study the decomposition of  the space ${\mathcal M}_n(\R^m, \C_m)$ 
of spherical monogenics of order $n$ under the action of $\Spin(p) \times \Spin(q)$. As a result we obtain a  $\Spin(p) \times \Spin(q)$-invariant orthonormal basis for ${\mathcal M}_n(\R^m, \C_m)$. 
In particular, using the construction with $p=2$ inductively, this yields a new orthonormal basis
for the space ${\mathcal M}_n(\R^m, \C_m).$



\end{abstract}
\maketitle

%%%%%%%%%%%%%%%%%%%%%%%%%%%%%%%%%%%%%%%%%%%%%%%%%%%%%%%%%%%%%%%%%%%%%%%%%%%%%%%%%%%%%%%%%%%%%%%%%%%%%%%%%%%%%%%%

\section{ \bf Introduction}

One of the main subjects studied in Clifford analysis is the function theory of monogenic functions and its interaction with the representation theory of the group $\Spin(m)$ (see e.g. \cite{BDS},\cite{DSS},\cite{GM}). 
Let  $ (e_1, \ldots , e_m) $ be an orthonormal basis of the Euclidean space $\R^m$. With respect to this basis, the Dirac operator on $\R^m$ is given by $ \px = \sum_j e_j \partial_{x_j} $.  A monogenic function $f$ is a Clifford algebra (or spinor)-valued function satisfying $\partial_{x}f=0$ in some open set $\Om\sub \R^m$. 

A crucial result is the fact that all possible (half-)integer irreducible representations of $\Spin(m)$ can be realized by means of so-called harmonic (monogenic) polynomials of several vector variables (see e.g. \cite{VLSC} and in particular for the harmonic case \cite{Vergne} or \cite{Goodman} ).  
In this paper we will focus our attention to the representations of highest weight $(n,0,\ldots, 0)$ and
$(n+\frac{1}{2},\frac{1}{2},\ldots, \frac{1}{2})$.  %; the so-called type I-representations. \\

The space $\Har_n(\R^m)$ of $n$-homogeneous harmonic polynomials in $\R^m$ is a well known model for the first type 
of representation. These are the well known spherical harmonics which play an important role in the harmonic analysis of the Laplace operator in $\R^m$ (see for example \cite{Müller}, \cite{Vilenkin1}).
The action of $s\in  \SO(m)$ on $P(x) \in \Har_n(\R^m)$ is given by $h(s)P(x)=P(s^{-1}xs)$ and the corresponding Lie algebra $\mathfrak{so}(m)$ is generated by the operators $L_{ij}=x_i\partial_{x_j}-x_j\partial_{x_i}$, $i<j$; $i,j=1, \ldots,m$. Put $N=\left[\frac{m}{2}\right]$. Fix the Cartan subalgebra $\mathfrak{h}=\{L_{2j-1,\,2j}:j=1,\ldots, N\}$ of $\mathfrak{so}(m)$. 

Using tools of Clifford analysis, the half-integer-representation can be described by means of spherical monogenics of order $n$. These are homogeneous $\R_m$ ($\C_m$)-valued monogenic polynomials of degree $n$ in $\R^m$. (If one wants irreducibility one needs to consider spinor-valued monogenics.) 
The corresponding spaces are denoted  by ${\mathcal M}_n(\R^m, \R_m)$ or ${\mathcal M}_n(\R^m, \C_m)$. They are a refinement of the notion of spherical harmonic in the sense that 
$$
\Har_n(\R^m, \C_m)={\mathcal M}_n(\R^m, \C_m) \oplus x{\mathcal M}_{n-1}(\R^m, \C_m)
$$ 
where left multiplication with the vector variable $x$ is to be considered as a $\Spin(m)$-invariant embedding. 
The action of $s\in  \Spin(m)$ on $P(x) \in \Mon_n(\R^m)$ is now given by a different representation $L(s)P(x)=sP(s^{-1}xs)$. The corresponding Lie algebra $\mathfrak{spin}(m)$ is generated by the momentum operators $M_{ij}=L_{ij}-\frac{1}{2}e_{ij}$, $i<j$; $i,j=1, \ldots,m$ where $e_{ij}=e_i \wedge e_j$. 
Here we choose the Cartan subalgebra $\mathfrak{h}=\{M_{2j-1,\,2j}:j=1,\ldots, N\}$ of $\mathfrak{spin}(m)$. 
Of course $\mathfrak{so}(m) \cong \mathfrak{spin}(m) \cong \mathfrak{o}(m)$ as Lie algebras. To stress the fact that our realizations of these Lie algebras arise from different representations we will use both notations $\mathfrak{so}(m)$ and $\mathfrak{spin}(m)$.

Fix the direct sum $\R^m=\R^p \oplus \R^q$. The aim of this paper is to study the decomposition of $\Har_n(\R^m)$ and in particular ${\mathcal M}_n(\R^m, \C_m)$ under the action of $\Spin(p) \times \Spin(q)$. We will prove that this decomposition (in a certain sense) is multiplicity free. Moreover we provide a construction of an orthonormal basis of $\Har_n(\R^m)$ and ${\mathcal M}_n(\R^m, \C_m)$. Another  construction of a $\Spin(p) \times \Spin(q)$-invariant basis of 
${\mathcal M}_n(\R^m, \C_m)$ was given in \cite{DSS} and in the paper \cite{So4} by F.~Sommen. 
The key ingredients in our construction are the harmonic and monogenic Fischer decomposition which can also be reformulated in the language  of the Howe dual pairs $(\SO(m),\mathfrak{sl}(2))$ and  $(\Spin(m),\mathfrak{osp}(1|2))$. Another basic tool is the extremal projector (see \cite{AST}, \cite{Zh1}) corresponding to the dual partners $\mathfrak{sl}(2)$ and $\mathfrak{osp}(1|2)$ of our initial groups $\SO(m)$ and $\Spin(m)$. Using the explicit version of this projection we obtain formulae for the basis vectors in terms of Jacobi polynomials. Similar expressions for the basis vectors were obtained in \cite{So4} by solving a system of differential equations. 
As special cases we will treat the cases $p=1,2$ in more detail.    

Let us consider first of all the harmonic case.
Let $p=1$. Take the standard ONB $(e_1, \ldots,e_m)$ of $\R^m$ and the corresponding chain of subgroups
\bear
\SO(m) \supset \SO(m-1) \supset \ldots \supset \SO(2) \,,
\enar
where $\SO(m-i)$ is the subgroup of $\SO(m)$ fixing the vectors $e_1, \ldots,e_i$. 
Each inclusion $\SO(i) \supset \SO(i-1)$ gives rise to a branching of a (irreducible) $\SO(i)$-representation under the action of $\SO(i-1)$. This branching is multiplicity free and by induction we thus obtain an orthonormal basis for $\Har_n(\R^m)$. It is a standard result that this can be done for each irreducible representation of $\SO(m)$ and one obtains an orthogonal basis which is indexed by Gel'fand-Zetlin patterns (see \cite{Gelfand1}).   
Next, consider $p=2$. Take the Cartan basis $\mathfrak{h}=\{L_{12}, L_{34}, \ldots, L_{2N-1\, 2N}\}$ of $\mathfrak{so}(m)$ and the corresponding chain 
\bear
\SO(m) \supset \SO(m-2) \supset \ldots \supset \SO(2) \mbox{ or } \SO(1)
\enar
compatible with the maximal torus corresponding to the choice of $\mathfrak{h}$. Depending on the parity of $m$ this chain stops at $\SO(2)$ ($m$ even) or $\SO(1)=\{1\}$ ($m$ odd).  
By induction we thus obtain an orthogonal basis of eigenfunctions of $\mathfrak{h}$ for $\Har_n(\R^m)$.  

Consider now the monogenic case and $p=1$.  Let $\R^m=\R\, e_1 \oplus e_1^{\perp} \cong \R \, e_1 \oplus \R^{m-1}$. The space of spherical monogenics 
$\Mon_k(\R, \C_1)$ on the line is non-trivial if $k=0$ and  $\Mon_0(\R, \C_1) \cong \C_1 \cong a+be_1,\, a,b \in \C$.
Consider the chain of subgroups 
\bear
\Spin(m) \supset \Spin(m-1) \supset \ldots \supset \Spin(2) 
\enar
in the maximal flag of subspaces 
\bear
\R^m  \supset e_1^{\perp} \supset \ldots \supset (e_{1} \wedge \ldots \wedge e_{m-2})^{\perp} \,. 
\enar
Similar to the harmonic case we obtain again a Gel'fand-Zetlin type basis for $\Mon_n(\R^m,\C_m)$. 
Next, take $p=2$. Let $\R^2=\mbox{span} \{e_1, e_2\}$. Then $\R^m=\R^2  \oplus \R^{m-2}$.
Take the Cartan basis $\mathfrak{h}=\{M_{12}, M_{34}, \ldots, M_{2N-1\, 2N}\}$ of $\mathfrak{spin}(m)$ and the corresponding chain 
\bear
\Spin(m) \supset \Spin(m-2) \supset \ldots \supset \Spin(2) \mbox{ or } \Spin(1)
\enar
in the maximal flag (preserved by the maximal torus)
\bear
\R^m  \supset (e_{1} \wedge e_2) ^{\perp} \supset \ldots \supset (e_{1} \wedge \ldots \wedge e_{m-2})^{\perp} \mbox{ or }  (e_{1} \wedge \ldots \wedge e_{m-1})^{\perp}\,. 
\enar
In this case, we reduce the symmetry from $\Spin(m)$ to $\Spin(2)\times \Spin(m-2).$ The branching rules for the product have the multiplicity one property, hence we can construct by induction
an analogue of the Gel'fand-Zetlin basis for the chain of inclusions
$$
\Spin(m)\supset\Spin(2)\times\Spin(m-2)\supset\Spin(2)\times\Spin(2)\times\Spin(m-4)\ldots
$$
This gives a new orthogonal basis for the space of spherical monogenics $\Mon_n(\R^m,\C_m)$.

The case $p=1$ (which we refer to as step one branching) and $p=2$ (step two branching) are quite different in the sense that they behave differently with respect to the Cartan basis $\mathfrak{h}$. 
If $p=1$ (the Gel'fand-Zetlin basis), the basis vectors are not eigenvectors
for the chosen Cartan subalgebras. This is due to the fact that the chain of $\SO$- or $\Spin$-groups 
alternates between two different types ($B$ and $D$) and the embeddings are not compatible with
the root systems. If $p=2$, the chains consist of Lie groups of the same type. 
The embeddings here respect the root systems so that the basis of $\Har_n(\R^m)$ and $\Mon_n(\R^m,\C_m)$ 
consists of weight vectors. We will show that for the type of representations of $\Spin(m)$ we consider in this paper, the basis which is induced by the procedure of step two branching is also orthogonal.

%%%%%%%%%%%%%%%%%%%%%%%%%%%%%%%%%%%%%%%%%%%%%%%%%%%%%%%%%%%%%%%%%%%%%%%%%%%%%%%%%%%%%%%%%%%%%%%%%%%%%%%%%%%%%%%%%%%%%%%%%%%%%

\section{ \bf Basic definitions from Clifford Analysis} 

In this section we collect some basic material and fix some notations. More detailed information  concerning Clifford algebras and Clifford analysis can be found in  
\cite{BDS},\cite{DSS},\cite{GM},\cite{guerleb}. 

Let $ (e_1, \ldots , e_m) $ be an orthonormal basis of Euclidean space 
$ \R^m $ endowed with the inner product $ \langle x,y  \rangle  
= \sum^m_{i=1} x_i y_i, \; x, y \in \R^m $. 
By $ \R_{0,m}$ we denote the real  $ 2^m$-dimensional 
Clifford algebra over $ \R^m$ generated by the relations 
$$ 
e_i e_j + e_j e_i = -2\delta_{ij}\,.
$$ 
If there is no confusion possible we use the shorter notation $ \R_{m}:= \R_{0,m}$. 
An element of $ \R_m $ is of the form  $ a = \sum_{A\subset M} a_A e_A, \; a_A \in \R, 
\; M = \{ 1, \ldots , m \} $ where $A$ is an ordered subset of $M$ and $ e_\phi = e_0 = 1 $. The $k$-vector part of $a$ is given by  $[a]_k=\sum_{|A|=k}a_A e_A$ and  $a=\sum_{k=0}^m [a ]_k$ with $ [a ]_k \in  \R_m^{(k)}$.  
Vectors $ x \in \R^m $ are identified with 1-vectors 
$ x = \sum^m_{j = 1} x_j e_j \in \R_m^{(1)}$. The Clifford product of two vectors $x$ and $y$ splits into minus the inner product and the wedge product of $x$ and $y$:
$$
xy=-\langle x,y \rangle + x \wedge y\,.
$$
The complex Clifford algebra $\C_m$ is the complexification of $\R_m$. 
The following (anti-)involutions are of importance. They are
defined by their action on the basis elements $e_i$ and extended by linearity to  $\R_m$:
\begin{itemize}
\item  main involution $a \mapsto a'$; \quad 
$
(ac)'=a'c'\,,\quad e'_i=-e_i\,,
$
\item  reversion $a \mapsto \tilde{a}$; \quad 
$
\quad \widetilde{ac}=\tilde{c}\tilde{a}\,,\quad \tilde{e}_i=e_i\,,
$
\item  conjugation $a \mapsto \bar{a}$; \quad
$
\quad \overline{ac}=\bar{c} \bar{a}\,,\quad \bar{e}_i=-e_i\,.
$
\end{itemize}
The main involution $'$ defines a $\Z_2$-grading on $\R_m$. The eigenspaces  $\R_m^{\pm}$ of the main involution $'$ are the so-called 
even and odd part of the Clifford algebra. The even part $\R_m^+$ is a subalgebra of $\R_m$ isomorphic to the Clifford algebra $\R_{m-1}$ and  $\R_m=\R_m^+ \oplus \R_m^-$.  
Conjugation  
on $ \C_m $ is the anti-involution on $ \C_m $ given by 
$ \bar{a} = 
\sum_{A \subset M} \bar{a}_A \bar{e}_A $.   

The following subgroups of the real Clifford algebra $\R_m$ are of
interest. The Pin group $\Pin(m)$ is the group consisting of products of 
unit vectors in $\R^m$; the Spin group $\Spin(m)$ is the subgroup of
$\Pin(m)$ consisting of products of an even number of unit vectors in
$\R^m$. For an element $s \in \Pin(m)$ the map $\chi(s):\R^m \ra \R^m:x
\mapsto sx (s')^{-1}$ induces an orthogonal transformation of $\R^m$. In
this way $\Pin(m)$ defines a double covering of the orthogonal group
$\mbox{O}(m)$. The restriction of this map to $\Spin(m)$
defines  
a double covering of the special orthogonal group $\SO(m)$.  

The Dirac operator on $\R^m$ is given by $ \px = \sum_{j=1}^m e_j \partial_{x_j} $ and has the fundamental property that $\px^2=-\lap_x$. Let $\Om$ be an open subset of $\R^m$ and let $f$ be a Clifford algebra (or spinor-valued) function; $f$ is said to be monogenic in $\Om$ if $\partial_x f=0$ in $\Om$ .    
The unit sphere in $ \R^m $ is denoted by $S^{m-1}$.  
Consider  polar coordinates $(\rho, \omega) \in \R_+ \times S^{m-1}$ in $\R^m$:
$$ 
x = \rho \omega , \; \rho = |x| = 
(x^2_1 + \ldots + x^2_m )^{1/2}\,, \quad \omega \in S^{m-1}\,. 
$$
The Dirac operator admits the polar 
decomposition 
$$ 
\px = \omega ( \partial_\rho + \frac{1}{\rho} 
\Gamma_{\omega} ) \quad  \mbox{with}  \quad \Gamma_{\omega} = - x \wedge \px  
$$
being the spherical 
Dirac operator on $ S^{m-1} $. 
In terms of the momentum operators
$L_{ij}=x_i \ran_{x_j}-x_j \ran_{x_j}$, the
$\G$-operator can be expressed as  $\G=-\sum_{i<j} e_{ij} L_{ij}$.   

Let $ n \in \N $. The space of $V$-valued $n$-homogeneous monogenic polynomials on $\R^m$ is denoted by ${\mathcal M}_n(\R^m, V)$. These type of polynomials are known as spherical monogenics of order $n$. For our purposes, $V$ is usually $\R_m$, $\C_m$ or a spinor space $\S$. Spherical monogenics can alternatively be defined by the conditions: 
$$
\Gamma P_n=-n P_n\,, \quad E P_n=nP_n\,; \quad E:=\inpro{x}{\partial_x}
$$ 
being the Euler operator on $\R^m$. Further information on spherical monogenics can be found in e.g. \cite{BDS},\cite{DSS},\cite{guerleb},\cite{So2},\cite{VL1} . 

The space of $\C$-valued $n$-homogeneous harmonic polynomials ($\lap_x P(x)=0$) on $\R^m$ is denoted by $\Har_n(\R^m, \C)$. Using the polar decomposition of the Laplacian
\bear
\lap_x=\partial_{\rho}^2+\frac{m-1}{\rho}\partial_{\rho}+\frac{1}{\rho^2}\lap_{LB} \, , 
\enar
spherical harmonics can alternatively be defined by the conditions: 
$$
\lap_{LB} P_n=-n(n+m-2) P_n\,, \quad E P_n=n P_n
$$  
where $\lap_{LB}$ has a two-fold meaning: it is either the Laplace-Beltrami operator on the sphere or the Casimir operator of the $h$-representation (see also next section).
%%%%%%%%%%%%%%%%%%%%%%%%%%%%%%%%%%%%%%%%%%%%%%%%%%%%%%%%%%%%%%%%%%%%%%%%%%%%%%%%%%%%%%%%%%%%%%%%%%%%%%%%%%%%%%%%%%%%%%%%%%%%%%%%%%%%%%%%%

\section{ \bf Representations of $\Spin(m)$}

The Lie algebra  $\mathfrak{spin}(m)$ can be realized  inside the space $\R_{m}^{(2)}$ of bivectors in $\R_m$ endowed with the usual commutator product $[\,,\,]$: 
$\mathfrak{spin}(m) \cong \Span \{-\frac{1}{2}e_{ij}: i<j, i,j=1,\ldots,m \}$. The bivectors $e_{ij}=e_i \wedge e_j$  satisfy the commutation relations 
\bear
\left[e_{ij}, e_{kl}\right]=2\delta_{il} e_{kj}+2\delta_{jl} e_{ik}-2\delta_{ik} e_{lj}-2\delta_{jk} e_{il}\,.
\enar 
Putting $B_{ij}=-\frac{1}{2}e_{ij}$ one obtains the usual defining relations of $\mathfrak{so}(m)$ (or here also frequently denoted as $\mathfrak{spin}(m)$):
\bear
\left[B_{ij}, B_{kl}\right]=-\delta_{il} B_{kj}-\delta_{jl} B_{ik}+\delta_{ik} B_{lj}+\delta_{jk} B_{il}\,.
\enar 
Let $\rho: \mbox{\Spin}(m) \ra \Aut (V)$ be a representation of $\mbox{\Spin}(m)$ on some vector space $V$.
The infinitesimal (derived) representation of $\rho$ on $V$ is given by
$$
d\rho(w)f=\lim_{t \ra 0}\frac{1}{t}(\rho(\exp(t w))-1)f\,, \quad  w \in \R_{m}^{(2)}, \, f \in V .
$$ 
The Casimir operator of the representation $\rho$ of $\Spin(m)$  is then defined by
$$
C(\rho)=\frac{1}{4}\sum_{i<j}d\rho(e_{ij})^2 .
$$
The Casimir operator $C(\rho)$ acts as a multiple of the identity on each
$\Spin(m)$-irreducible piece occurring in $V$. 
Take now for $V$  the space of $\C_m$-valued polynomials ${\mathcal P}(\R^m, \C_m)$ and take $s \in \Spin(m)$. 
Consider the 
following two unitary (for the Fischer inner product) representations of $\Spin(m)$ on ${\mathcal P}(\R^m, \C_m)$:
\bear
H(s)P(x)&=& s P(s^{-1}x s)s^{-1}=s P(\bar{s}x s)\bar{s}
\\ 
L(s)P(x)&=& s P(s^{-1}x s)=s P(\bar{s}x s)\,.
\enar
The $L$-representation typically acts on (spinor-valued) monogenic polynomials.
The action $a \mapsto s a s^{-1}$ preserves the space $\C_m^{(k)}$ of $k$-vectors, hence the $H$-representation may also act on the space 
${\mathcal P}(\R^m, \C_{m}^{(k)})$ of $\C_{m}^{(k)}$-valued 
polynomials. This defines for each $k=0, \ldots, m$ a representation of $\SO(m)$. 
%For our purposes, the $L$-representation will act on monogenic polynomials.  
In case $H$ acts on scalar-valued polynomials ($k=0$) one also uses the notation $h$ instead of $H$. We thus get the usual representation $h$ of $\SO(m)$ on ${\mathcal P}(\R^m, \C)$ and 
\bear
dh(e_{ij})=-2L_{ij}:=-2(x_i \partial_{x_j}-x_j \partial_{x_i})\,, \quad  \quad dL(e_{ij})=-2L_{ij}+e_{ij}:=-2M_{ij}\,.
\enar
The corresponding Casimir operators are 
\bear
C(h)&=&\sum_{i<j}L^2_{ij}=\lap_{LB}  \\
C(L)&=&\sum_{i<j}M^2_{ij}=\lap_{LB} +\Gamma -\frac{1}{4} \bin{m}{2} 
                       =\Gamma(m-1-\Gamma) -\frac{1}{4} \bin{m}{2}
\enar
%where $L_{ij}$, $M_{ij}$ are the (spin)-momentum operators and $\lap_{LB}$ is the Laplace-Beltrami operator on $\sfe$.      
The Cartan subalgebra $\mathfrak{h}$ of $\mathfrak{spin}(m)$ is fixed by 
$\mathfrak{h}=\{M_{2j-1,\,2j}:j=1,\ldots, N=\left[\frac{m}{2}\right]\}$.  The exponential of $\mathfrak{h}$  yields the corresponding maximal torus
\bear
H=\{\exp(-\frac{1}{2}e_{12} t_1) \ldots \exp(-\frac{1}{2}e_{2N-1\,2N}t_N), t_i \in [0, 2\pi [ \}, \quad N=\left[\frac{m}{2}\right]\,.
\enar 

\section{ \bf Spherical monogenics and representations}

Models for irreducible representations of $\Spin(m)$ can be realized by means of monogenic polynomials of several vector variables (see \cite{VLSC}). 
We will illustrate this for the type of representations which appear in this paper. First of all, models for the spinor representation can be realized inside the complex Clifford algebra $\C_m$. The action $L(s) a = s a$ of $s \in \Spin(m)$ on $\C_m$ 
leads to the fundamental representation  of $\Spin(m)$ on the spinor space
$\S$. A model for $\S$ can be constructed as follows.  
Let $m=2N+1$ be odd. Consider the Witt basis
$$                                       
T_{j,+} := \frac{1}{2}(e_{2j-1} + i e_{2j})\,,\quad T_{j, -} := \frac{1}{2}(e_{2j-1} - i e_{2j}) 
$$                             
and the idempotents  $I_{j, +} = -T_{j,+}T_{j,-}$. The product 
$I_+ = I_{1, +} \ldots I_{N, +}$ 
defines a  primitive idempotent; the ideal $\C^+_m I_+$ is minimal and gives a
model for the spinor space $\S$. 
The action of the maximal torus gives
\begin{eqnarray*}
L(s) I_+ &=& \exp (-\frac{1}{2} (t_1 e_{12} + \ldots + t_N e_{2N-1, 2N})) I_+\\
&=& \exp (\frac{i}{2}(t_1 + \ldots + t_N)) I_+;
\end{eqnarray*}
and the weight is given by $(\frac{1}{2}, \ldots,
\frac{1}{2})$.
A model for the irreducible representation with weight $(n+\frac{1}{2},\frac{1}{2}, \ldots, \frac{1}{2})$ can be obtained by considering the $L$-action of $\Spin(m)$ on the highest weight vector
$$
(x_1+i x_2)^n I_+ \,.
$$
For $s$ of the maximal torus $H$, we have that
\begin{eqnarray*}
L(s) (x_1+i x_2)^n I_+ &=& \exp (-n t_1 e_{12}) \exp (-\frac{1}{2} (t_1 e_{12} + \ldots + t_N e_{2N-1, 2N}))(x_1+i x_2)^n  I_+\\
&=& \exp (in t_1)\exp (\frac{i}{2}(t_1 + \ldots + t_N)) (x_1+i x_2)^n  I_+ \,.
\end{eqnarray*}
The resulting irreducible $\Spin(m)$-module is the space of spinor-valued spherical monogenics of order $n$:
$$
{\mathcal M}_n(\R^m, \S):={\mathcal M}_n(\R^m, \C_m^+)I_+  \cong (n+\frac{1}{2},\frac{1}{2}, \ldots, \frac{1}{2})\,.
$$ 
Thus ${\mathcal M}_n(\R^m, \C_m)$ is the direct sum of several equivalent representations with weight 
$(n+\frac{1}{2},\frac{1}{2},\ldots, \frac{1}{2})$.   \\
The even dimensional case $(m=2N)$ requires a little bit more care because there are now two inequivalent 
spinor representations of $\Spin(m)$. 
Consider now the spinor spaces $\S_+ \cong \C^+_m I_+$ and $\S_- \cong \C^+_m I_-$ where the primitive 
idempotents  $I_+$ and $I_-$
are given by
$$                                                       
I_+ := I_{1,+} \ldots I_{N-1, +} I_{N, +}\,, \quad I_- := I_{1,+} \ldots I_{N-1, +} I_{N, -}\,,
$$
and 
$$
I_{N,-} := \frac{1}{2}(1 - i e_{m-1} e_m). 
$$
As the pseudoscalar $e_M:= e_1 \ldots e_m$ is 
actually $\Spin(m)$-invariant (but not $\Pin(m)$-invariant) and $e_M^2=(-1)^N$, there are 
two $\Spin(m)$-invariant projections
$$
P_+ = \frac{1}{2}(1 + i^N e_M) \quad \mbox{and} \quad P_- = \frac{1}{2}(1 - i^N e_M)
$$
onto the $\pm$-eigenspaces of $i^N e_M$. Now $I_+ = P_+ I_+$ and $I_- = P_- I_-$, hence 
the spinor spaces $\C^+_m I_+$ and $\C^+_m I_-$ are inequivalent under the
action  of $\Spin(m)$. The weights are obtained
from the action of the maximal torus and are given by 
$(\frac{1}{2}, \ldots,\frac{1}{2}, \frac{1}{2}) $ resp. $(\frac{1}{2}, \ldots,\frac{1}{2},
-\frac{1}{2})$. 

A model for the irreducible representation with weight $(n+\frac{1}{2},\frac{1}{2}, \ldots, \pm \frac{1}{2})$ can be obtained by considering the $L$-action of $\Spin(m)$ on the highest weight vectors
$$
(x_1+i x_2)^n I_+ \quad \mbox{ and } \quad (x_1-i x_2)^n I_- \,.
$$
In the first case, the action of the maximal torus yields the same weight as in the odd dimensional case. In the second case, the last entry in the weight changes sign because now
\begin{eqnarray*}
L(s) (x_1-i x_2)^n I_- &=& \exp (-n t_1 e_{12}) \exp (-\frac{1}{2} (t_1 e_{12} + \ldots + t_N e_{2N-1, 2N}))(x_1-i x_2)^n I_-\\
&=& \exp (in t_1)\exp (\frac{i}{2}(t_1 + \ldots + t_{N-1}- t_N)) (x_1-i x_2)^n  I_- \,.
\end{eqnarray*}
The resulting non-equivalent irreducible $\Spin(m)$-modules are the spaces of spinor ($\S^+$ or $\S^-$)-valued spherical monogenics of order $n$:
$$
{\mathcal M}_n(\R^m, \S^+):={\mathcal M}_n(\R^m, \C_m^+)I_+ \quad \mbox{ and } \quad {\mathcal M}_n(\R^m, \S^-):={\mathcal M}_n(\R^m, \C_m^+)I_-\,.
$$

%%%%%%%%%%%%%%%%%%%%%%%%%%%%%%%%%%%%%%%%%%%%%%%%%%%%%%%%%%%%%%%%%%%%%%%%%%%%%%%%%%%%%%%%%%%%%%%%%%%%%%%%%%%%%%%%%%%%%%%%%%%%%%%%%%%%%%%%%%%%%%%%%%%%%%%%%%%%%%% 

\section{\bf Inner products on  ${\mathcal M}_n(\R^{m},\C_{m})$ }
In the context of spherical monogenics one usually considers the so-called Fischer inner product or the standard inner product on the sphere. 
Let $P,Q$ be $\C_m$-valued polynomials on $\R^m$. The Fischer inner product of $P$ and $Q$ is the positive definite Hermitean inner product defined as
$$
\inpro{P}{Q}_m:=[\bar{P}(\partial_x)Q(x)]_0|_{x=0}   \,.
$$
This definition implies immediately that polynomials of different degree are Fischer orthogonal. 
The adjoint of $A \in  \End({\mathcal P}(\R^{m},\C_{m}))$ relative to the Fischer inner product is denoted by $A^{\adj}$.
Up to a sign the variables $x_i$ and the derivatives $\partial_{x_i}$ are Fischer-adjoint while the generators $e_i$ of the Clifford algebra $\C_m$ are skew-adjoint:
\bear
\inpro{x_iP}{Q}_m  &=& -\inpro{P}{\partial_{x_i}Q}_m \\
\inpro{e_i P}{Q}_m  &=& -\inpro{P}{e_i Q}_m \,.
\enar
Therefore the Fischer-adjoint of $|x|^2$ is the Laplace operator $\lap_x$ on $\R^m$ and the Fischer-adjoint of the vector variable $x$ is the Dirac operator $\partial_x$ on $\R^m$. This plays a crucial role in the Fischer decomposition which will be explained later on. 
Call $A_m$  the surface area of the unit sphere $S^{m-1}$. 
The normalized $\C_m$-valued inner product on $L_2(S^{m-1},\,{\mathbb C}_{m})$ is given by:
\bear
\inpro{f}{g}=\frac{1}{A_m}\int_{S^{m-1}}\bar{f}(\omega)g(\omega)\,dS(\omega)\,.
\enar
This inner product can be turned into a $\C$-valued inner product by projecting it onto its scalar part. The corresponding inner product will be denoted as $\inpro{\,}{\,}_{S^{m-1}}$. All of these inner products are $\Spin(m)$-invariant. As the space  ${\mathcal M}_n(\R^{m},\C_{m})$ is basically irreducible under the action of $\Spin(m)$, it is not surprising that 
the Fischer inner product on ${\mathcal M}_n(\R^{m}, \C_m)$ and the inner product on the sphere are proportional. They are related as follows:
\bear
\inpro{P}{Q}_m=\frac{2^n \Gamma(n+\frac{m}{2})}{\Gamma(\frac{m}{2})} \inpro{P}{Q}_{S^{m-1}}\,,\quad P,Q \in {\mathcal M}_n(\R^{m}, \C_m)\,.
\enar
More details of the notions presented here can also be found in e.g. \cite{GM},\cite{DSS2},\cite{DSS},\cite{guerleb}.

%%%%%%%%%%%%%%%%%%%%%%%%%%%%%%%%%%%%%%%%%%%%%%%%%%%%%%%%%%%%%%%%%%%%%%%%%%%%%%%%%%%%%%%%%%%%%%%%%%%%%%%%%%%%%%%%%%%%%%%%%%%%%%%%%%%%%%%%%

\section{ \bf The harmonic projection}
The following Lemma lists the (anti-)commutation relations between various $\Spin(m)$-invariant operators which play a fundamental role in classical Clifford analysis in one vector variable. They follow from standard calculations (see also \cite{BSES}). 
\begin{lemma} $($Basic $($anti-$)$commutation relations in $\R^m)$ 

\bear
\begin{array}{ll}
\left[\lap_x, |x|^2 \right]=4(E_x+\frac{m}{2}) & \hspace*{10mm}\{x, \partial_x\}=-m-2E_x \\[2mm]
\left[E_x+\frac{m}{2}, |x|^2\right]=2|x|^2 & \hspace*{10mm} \left[E_x+\frac{m}{2}, x\right]=x \\[2mm]
\left[E_x+\frac{m}{2}, \lap_x\right]=-2\lap_x  & \hspace*{10mm}\left[E_x+\frac{m}{2}, \partial_x\right]=-\partial_x \\[2mm]
\left[\partial_x, |x|^{2l}\right]=2l\,|x|^{2l-2} x & \hspace*{10mm}\left\{\partial_x, |x|^{2l}\right\}=(-m-2E_x+2l)|x|^{2l-2} \\[2mm]
\left[\lap_x, |x|^{2j} \right]=4j(E_x+\frac{m}{2}-j+1)|x|^{2j-2} %\hspace*{10mm} \phantom{\left\{\partial_u, |u|^{2l}\right\}} 
\end{array}
\enar
\end{lemma}
Identify 
\bea
\Bp:=-\frac{1}{2}\lap_x\,, \quad \Bm:=\frac{1}{2}|x|^2\,, \quad \Bn:=-(E_x+\frac{m}{2})\,, \quad \Fp:=-\frac{1}{\sqrt{2}}\partial_x\,, \quad \Fm:=-\frac{1}{\sqrt{2}}x \label{identification}\,.
\ena
Then $\{\Bp, \Bm, \Bn, \Fp, \Fm \}$ satisfy the (anti-)commutation relations of the Lie super algebra $\gosp(1|2)$. As a Lie super algebra this algebra is generated by the odd (fermionic) generators $\Fp, \Fm $. The even (bosonic) elements $\Bp, \Bm, \Bn$ generate the even Lie subalgebra $\mathfrak{sl}(2)$. 
Let us first consider the Fischer decomposition related to the Laplace operator: this means that we are only considering the even part $\mathfrak{sl}(2)$ of $\gosp(1|2)$. We have
\bear
\Pol(\R^m)&=& \bigoplus_{s=0}^{\infty} |x|^{2s} \Har(\R^m) \\
          &=& \bigoplus_{k=0}^{\infty} \bigoplus_{s=0}^{\infty} |x|^{2s} \Har_k(\R^m) \\
          & \cong & \bigoplus_{k=0}^{\infty} \Ve_k \otimes \Har_k(\R^m) \,.
\enar 
Here $\Ve_k$ is the irreducible $\mathfrak{sl}(2)$-module with weight $-(k+\frac{m}{2})$. 
The last isomorphism indicates that the space $\Pol(\R^m)$ has a multiplicity free decomposition under the joint action $\mathfrak{sl}(2) \times \SO(m)$. The pair $(\mathfrak{sl}(2), \SO(m))$ is a particular example of a Howe dual pair. 
Each $H(x) \in \Har(\R^m)$ can be regarded as a highest weight vector annihilated by the positive root $\Bp:=-\frac{1}{2}\lap_x$.  The corresponding $\mathfrak{sl}(2)$-module generated by  $H_k(x) \in \Har_k(\R^m)$ is the infinite dimensional module  $\bigoplus_{s=0}^{\infty} |x|^{2s} H_k(x)$. 
Given a polynomial $P \in \Pol(\R^m)$ we thus have a unique decomposition
\bear
P(x)=\sum_{j=0}^{\infty} |x|^{2j} H^{(j)}(x)
\enar
where $H^{(j)}(x) \in \Har(\R^m)$. The harmonic polynomial $H^{(0)}(x)$ is called the harmonic part of $P(x)$ and the corresponding map 
$\phar: \Pol(\R^m) \ra \Har(\R^m)$ is the harmonic projection as mentioned in  \cite{Vilenkin1}. In the same way one can define 
the projections $P_{H,2s}:\Pol(\R^m) \mapsto |x|^{2s}\Har(\R^m)$. These projections are closely related to the harmonic projection $\phar$. Obviously $P_{H,0}=\phar$.  We will now give a formula to determine $P_{H,2s}$. Because of its importance we derive the formula for our specific situation. 
We also would like to point out that the explicit expression for $\phar$ is not new and appears at various places in the literature. 

Another useful interpretation is that $P_H$ is a specific realization of the extremal projector of the Lie algebra $\mathfrak{sl}(2)$. 
Consider the Gauss decomposition $\mathfrak{sl}(2)=\mathfrak{n}_- +\mathfrak{h}+\mathfrak{n}_+$. For an $\mathfrak{sl}(2)$-module $V$, the extremal projector $P$ is an operator which is constructed in a suitable extension of $U(\mathfrak{sl}(2))$ and projects $V$ onto its subspace $V^{\mathfrak{n}_+}$ (of highest weight vectors) parallel to $\mathfrak{n}_-V$. The extremal  projectors corresponding to the classical Lie algebras and more general algebraic structures have been studied by various authors (see \cite{Tol1}, \cite{Zh1}).

%%%%%%%%%%%%%%%%%%%%%%%%%%%%%%%%%%%%%%%%%%%%%%%%%%%%%%%%%%%%%%%%%%%%%%%%%%%%%%%%%%%%%%%%%%%%%%%%%%%%%%%%%%%%%%%%%%%%%%%%%%%%%%%%%%%%%%%%%%%%%%%%%%%%%%%%%%%%%%%%
\begin{theorem} 
Consider the Fischer orthogonal direct sum $\Pol(\R^m)= \bigoplus_{s=0}^{\infty} |x|^{2s} \Har(\R^m) $. 
The harmonic projection $\phar: \Pol(\R^m) \ra \Har(\R^m)$ and the projections $P_{H,2s}:\Pol(\R^m) \mapsto |x|^{2s}\Har(\R^m)$ can be expressed  by means of the  operators 
\bear
\phar&=&\sum_{j=0}^{\infty}  \frac{1}{2^{2j}j \,!}\frac{\G(-E_x-\frac{m}{2}+2)}{\G(-E_x-\frac{m}{2}+j+2)} |x|^{2j}\lap_x^j \\
     &=&\sum_{j=0}^{\infty}  \frac{1}{2^{2j}j \,!}\kappa_j(-E_x-\frac{m}{2}) |x|^{2j}\lap_x^j   \\
P_{H,2s}
   &=&A_s(E_x-2s)|x|^{2s} \phar \lap_x^s  \,,
\enar
where 
\bear
\kappa_j(z)&:=&\frac{\G(z+2)}{\G(z+j+2)}=\frac{1}{(z+2) \ldots(z+j+1)}=\frac{1}{(z+2)_j}  \\
   A_s(E_x)&:=& \frac{1}{2^{2s} s!\,} \frac{1}{(E_x+\frac{m}{2})_s} \,.
\enar
\end{theorem}
%%%%%%%%%%%%%%%%%%%%%%%%%%%%%%%%%%%%%%%%%%%%%%%%%%%%%%%%%%%%%%%%%%%%%%%%%%%%%%%%%%%%%%%%%%%%%%%%%%%%%%%%%%%%%%%%%%%%%%%%%%%%%%%%%%%%%%%%%%%%%%%%%%%%%%%%%%%%%%%%%%%%%
\bewijs
As ansatz we take 
\bea
P=\sum_{j=0}^{\infty} K_j(E_x) |x|^{2j}\lap_x^j \label{P-ansatz}
\ena
where $K_j(E_x)$ are unknown functions (not necessarily polynomials) of $E_x$. Thus $P$ rather belongs to some extension of $U(\mathfrak{sl}(2))$. We will now determine the solution of $\lap_x P=0$ in this extension (which will be described afterwards). Since
\bear
\lap_x P&=&\sum_{j=0}^{\infty} K_j(E_x+2) \lap_x |x|^{2j}\lap_x^j   \\
        %&=&\sum_{j=1}^{\infty} K_j(E_u+2) 4j(E_u+\frac{m}{2}-j+1)|u|^{2j-2}\lap_u^{j}+\sum_{j=0}^{\infty} K_{j}(E_u+2)|u|^{2j}\lap_u^{j+1}  \\
        &=&\sum_{j=0}^{\infty} \left(4(j+1)K_{j+1}(E_x+2)(E_x+\frac{m}{2}-j)+ K_{j}(E_x+2)\right) |x|^{2j}\lap_x^{j+1}   
\enar
one obtains that $\lap_x P=0$ if 
\bear
\frac{K_{j+1}(E_x+2)}{K_j(E_x+2)}=-\frac{1}{4(j+1)(E_x+\frac{m}{2}-j)}\,.
\enar
We thus obtain the unique solution (if we put $K_0=1$):
\bear
K_j(E_x+2)=\frac{(-1)^j}{2^{2j}j \,!} \frac{1}{(E_x+\frac{m}{2}-j+1) \ldots (E_x+\frac{m}{2})} \,.
\enar
The function $K_j(E_x)$ turns out to be  a rational function in $E_x$. The Poincar\'{e}-Birkhoff-Witt-Theorem tells us that $U(\mathfrak{sl}(2))$ as a vector space has a basis consisting of the elements $\Bn^{i} \Bm^j \Bp^k $. Denoting the space of rational functions in $E_x$ by $R(E_x) \cong \C(E_x)$ we see that the appropriate extension of $U(\mathfrak{sl}(2))$ is of the form $R(\mathfrak{h}) \otimes_{\mathfrak{h}} U(\mathfrak{sl}(2))$. This space has a basis of the form  $R(\Bn) \Bm^j \Bp^k $  
where $R(\Bn)$ is a rational function in $\Bn$; the expression (\ref{P-ansatz}) is precisely of this form. Hence 
\bear
K_j(E_x)%&=&\frac{(-1)^j}{2^{2j}j \,!} \frac{1}{(E_x+\frac{m}{2}-j-1) \ldots (E_x+\frac{m}{2}-2)}  \\
          &=&\frac{(-1)^j}{2^{2j}j \,!}\frac{1}{(E_x+\frac{m-2}{2}-j)_j}  \\
         % &=&\frac{(-1)^j}{2^{2j}j \,!}\frac{\G(E_x+\frac{m-2}{2}-j)}{\G(E_x+\frac{m-2}{2})}    \\
         % &=&\frac{1}{2^{2j}j \,!} \frac{1}{(-E_x-\frac{m}{2}+2) \ldots (-E_x-\frac{m}{2}+j+1)}  \\
          &=&\frac{1}{2^{2j}j \,!}\frac{1}{(-E_x-\frac{m}{2}+2)_j}  \\
          &=&\frac{1}{2^{2j}j \,!}\frac{\G(-E_x-\frac{m}{2}+2)}{\G(-E_x-\frac{m}{2}+j+2)} \,.
\enar
Here $\G(z)$ denotes the Gamma-function. This proves the statement.  

Next we determine $P_{H,2s}$. Remark that $\lap_x^s$ preserves the Fischer decomposition:
\bear
\lap_x^s: \sum_{i=0}^{\infty} |x|^{2i} H^{(i)}(x)  \ra \sum_{i=s}^{\infty} |x|^{2i-2s} H^{(i)}(x)\,, \quad  H^{(i)} \in \Har(\R^m) \,.
\enar
Hence $P_{H,2s}$ is up to some element of $\End(\Har(\R^m))$ of the form $|x|^{2s} \phar \lap_x^s$. Now take $H(x) \in \Har(\R^m)$. Then
\bear
\lap_x^s |x|^{2s}H(x)%&=& 2^{2s} s!\, (E_u+\frac{m}{2}) \ldots (E_u+\frac{m}{2}+s-1) H(u) \\
                     &=& 2^{2s} s!\, (E_x+\frac{m}{2})_s \, H(x)  
                     =\frac{1}{A_s(E_x)}  H(x)
\enar
if we put 
\bear
A_s(E_x):= \frac{1}{2^{2s} s!\,} \frac{1}{(E_x+\frac{m}{2})_s} \,.
\enar
The corresponding projection is now given by
\bear
P_{H,2s}&=&|x|^{2s} \phar A_s(E_x)\lap_x^s  \\
   &=&A_s(E_x-2s)|x|^{2s} \phar \lap_x^s  \,,
\enar
which completes the proof.
\eop

Recall the identification (\ref{identification}). Denoting
\bear
\kappa_j(\Bn)&=&\frac{1}{(\Bn+2) \ldots (\Bn+j+1)}  
       =\frac{\G(\Bn+2)}{\G(\Bn+j+2)}
\enar
we have 
\bear
\phar=\sum_{j=0}^{\infty}  \frac{(-1)^j}{j \,!}\kappa_j(\Bn) \Bm ^{j} \Bp ^j \,.
\enar
This form of $\phar$ is the so-called extremal projector for $\mathfrak{sl}(2)$ (see e.g. the work of Tolstoy \cite{Tol1}, Zhelobenko \cite{Zh1}). 
As mentioned before, the extremal projector $\phar$ does not belong to $U(\mathfrak{sl}(2))$ but it belongs to some
extension $TU(\mathfrak{sl}(2))$ of the universal enveloping algebra. A series $S$ belongs to the space $TU(\mathfrak{sl}(2))$ of  formal Taylor series 
if  
\bear
S=\sum_{i,j=0}^{\infty}R_{i,j}(\Bn) \Bm^i \Bp^j
\enar
where $R_{i,j}(\Bn)$ are rational functions of the Cartan element $\Bn$ and such that
there exists a natural number $n$ for which
$|i - j| \leq n$. 
One can show that $TU(\mathfrak{sl}(2))$ is an associative algebra with respect to the multiplication of formal series. 
The enveloping algebra $U(\mathfrak{sl}(2))$ has no zero divisors and therefore contains no non-trivial projections. As a result 
the equations
\bear
\Bp P = P \Bm = 0
\enar
have no non-trivial solutions in $U(\mathfrak{sl}(2))$. 
The algebra $TU(\mathfrak{sl}(2))$ on the other hand does contain non-trivial projections and the following system of equations 
\bear
\Bp P = P \Bm = 0\, , \quad \left[\Bn, P \right] = 0\, , \quad P^2 = P
\enar
has precisely the extremal projector $\phar$ as its  unique solution in $TU(\mathfrak{sl}(2))$. 
If we put 
$$
\alpha_s(\Bn)=(-1)^s 2^{2s}A_s(-\Bn-m-2s)
$$
we have under the usual identification
\bear
P_{H,2s}= \alpha_s(\Bn) \Bm^s \phar \Bp^s  \,.
\enar

\section{ \bf The monogenic projection}

We now consider the Fischer decomposition for the Dirac operator. Here we use the short notation $\Pol(\R^m)=\Pol(\R^m, \C_m)$, $\Mon_k(\R^m)=\Mon_k(\R^m, \C_m)$.  
We have:
\bear
\Pol(\R^m)&=& \bigoplus_{s=0}^{\infty} x^s \Mon(\R^m)  \\
          &=& \bigoplus_{j=0}^{\infty} |x|^{2j} \Mon(\R^m) \oplus \bigoplus_{j=0}^{\infty} x|x|^{2j} \Mon(\R^m) \\
          &=& \bigoplus_{k=0}^{\infty} \bigoplus_{s=0}^{\infty} x^{s} \Mon_k(\R^m) \\
          & \cong & \bigoplus_{k=0}^{\infty} \We_k \otimes \Mon_k(\R^m)  \,.
\enar 
This decomposition is a refinement of the harmonic (or $\mathfrak{sl}(2)$-) case in the sense that
\bear
\Har(\R^m)=\Mon(\R^m) \oplus x \Mon(\R^m)
\enar
Here $\We_k$ is an irreducible $\gosp(1|2)$-module with weight $-(k+\frac{m}{2})$. As an $\mathfrak{sl}(2)$-module $\We_k$ splits into two irreducible $\mathfrak{sl}(2)$-modules with weights $-(k+\frac{m}{2})$ and $-(k+1+\frac{m}{2})$. If we consider the space $\Pol(\R^m, \S)$ of spinor-valued polynomials, one obtains
\bear
\Pol(\R^m, \S) = \bigoplus_{k=0}^{\infty} \We_k \otimes \Mon_k(\R^m, \S)\,.
\enar 
This decomposition is once again multiplicity free if we now consider the joint action $\mathfrak{osp}(1|2) \times \Spin(m)$, thus providing another particular example of a Howe dual pair. Each $M_k(x) \in \Mon_k(\R^m)$ can be regarded as a highest weight vector annihilated by the positive root $\Fp:=-\frac{1}{\sqrt{2}}\partial_x$.  The corresponding $\gosp(1|2)$-module generated by  $M_k(x)$ is the infinite dimensional module  $\bigoplus_{s=0}^{\infty} x^{s} M_k(x)$ which contains the two inequivalent $\mathfrak{sl}(2)$-modules $\bigoplus_{s=0}^{\infty} x^{2s} M_k(x)$ and
$\bigoplus_{s=0}^{\infty} x^{2s} x M_k(x)$. Given a polynomial $P \in \Pol(\R^m)$ we thus have a unique decomposition
\bear
P(x)&=&\sum_{s=0}^{\infty} x^{s} M^{(s)}(x)   \\
&\cong &\sum_{j=0}^{\infty} |x|^{2j} M^{(2j)}(x) + \sum_{j=0}^{\infty} x|x|^{2j} M^{(2j+1)}(x) \\
\enar
where $M^{(s)}(x) \in \Mon(\R^m)$. The monogenic polynomial $M^{(0)}(x)$ is called the monogenic part of $P(x)$ and the corresponding map 
$\pmon:\Pol(\R^m) \ra \Mon(\R^m)$ is the monogenic projection. In the same way one can define two other type of projections 
\bear
P_{M,2s}&:&\Pol(\R^m) \ra |x|^{2s}\Mon(\R^m)  \\
P_{M,2s+1}&:&\Pol(\R^m) \ra |x|^{2s}x\Mon(\R^m)
\enar
These projections are clearly closely related to the monogenic projection $\pmon$. 
%%%%%%%%%%%%%%%%%%%%%%%%%%%%%%%%%%%%%%%%%%%%%%%%%%%%%%%%%%%%%%%%%%%%%%%%%%%%%%%%%%%%%%%%%%%%%%%%%%%%%%%%%%%%%%%%%%%%%%%%%%%%%%%%%%%%%%%%%%%%%%%%%%%%%%%%%%%%%%%%%%%
\begin{theorem} Consider the Fischer orthogonal direct sum 
\bear
\Pol(\R^m)=\bigoplus_{s=0}^{\infty} |x|^{2s} \Mon(\R^m) \oplus \bigoplus_{s=0}^{\infty} x|x|^{2s} \Mon(\R^m) 
\enar
The monogenic projection $\pmon:\Pol(\R^m) \ra \Mon(\R^m)$ is given by the operator
\bear
\pmon=\frac{x  \partial_x+2E_x+m-2}{2E_x+m-2}\phar \,
\enar
and the other projections can be expressed as 
\bear
P_{M,2s}
   &=&A_s(E_x-2s)|x|^{2s} \pmon \lap_x^s  \\ 
   %&=&A_s(E_u-2s)u^{2s} \pmon \partial_u^{2s}  \\
P_{M,2s+1}
   &=&-\frac{A_s(E_x-2s-1)}{2(E_x+\frac{m}{2}-s-1)}|x|^{2s} x \pmon \lap_x^s \partial_x 
   %&=&-\frac{A_s(E_u-2s-1)}{2(E_u+\frac{m}{2}-s-1)}u^{2s+1}  \pmon \partial_u^{2s+1}  \,.   
\enar
\end{theorem}
%%%%%%%%%%%%%%%%%%%%%%%%%%%%%%%%%%%%%%%%%%%%%%%%%%%%%%%%%%%%%%%%%%%%%%%%%%%%%%%%%%%%%%%%%%%%%%%%%%%%%%%%%%%%%%%%%%%%%%%%%%%%%%%%%%%%%%%%%%%%%%%%%%%%%%%%%%%%%%%%%%%%%%
\bewijs
Put $H(x)=\phar P(x)$, then $H(x)=M^{(0)}(x)+x M^{(1)}(x)$. Now $\partial_x H(x)= \partial_x x M^{(1)}(x)=-(m+2E_x)M^{(1)}(x)$. Hence
\bear
M^{(1)}(x)=-\frac{1}{m+2E_x}  \partial_x H(x)\,, \quad M^{(0)}=\frac{1}{m+2E_x-2}(m+2E_x-2+x  \partial_x) H(x)
\enar
and the monogenic projection takes the form
\bear
\pmon&=&\frac{1}{m+2E_x-2}(m+2E_x-2+x  \partial_x)\phar  \\
     &=&\phar +\frac{1}{m+2E_x-2}(x  \partial_x)\phar   \,.
\enar
To determine $P_{M,2s}$ remark that $\lap_x^s$ preserves the Fischer decomposition:
\bear
\lap_x^s: \sum_{i=0}^{\infty} x^{i} M^{(i)}(x)  \ra \sum_{i=2s}^{\infty} x^{i-2s} M^{(i)}(x)\,, \quad  M^{(i)}(x) \in \Mon(\R^m) \,.
\enar
Therefore $P_{M,2s}$ is up to some element of $\End(\Mon(\R^m))$ of the form $|x|^{2s} \pmon \lap_x^s$. Now take $M(x) \in \Mon(\R^m)$. Then
\bear
\lap_x^s |x|^{2s}M(x)%&=&2^{2s} s!\, (E_x+\frac{m}{2}) \ldots (E_x+\frac{m}{2}+s-1) M(x) \\
                     %=&2^{2s} s!\, (E_x+\frac{m}{2})_s \, M(x)  \\ 
                     =\frac{1}{A_s(E_x)}  M(x)  \,.
\enar
Hence
\bear
P_{M,2s}&=&|x|^{2s} \pmon A_s(E_x)\lap_x^s  \\
   &=&A_s(E_x-2s)|x|^{2s} \pmon \lap_x^s  
   %&=&A_s(E_x-2s)x^{2s} \pmon \partial_x^{2s}  \,.
\enar
Next, under the action of $\lap_x^s \partial_x $:
\bear
\lap_x^s \partial_x: \sum_{i=0}^{\infty} x^{i} M^{(i)}(x)  \ra \sum_{i=2s+1}^{\infty} x^{i-2s-1} M^{(i)}(x)\,, \quad M^{(i)}(x) \in \Mon(\R^m)  \,.
\enar
Now $P_{M,2s+1}$ is up to some element of $\End(\Mon(\R^m))$ of the form $|x|^{2s}x \pmon \lap_x^s \partial_x$ and:
\bear
\lap_x^s \partial_x x |x|^{2s}M(x)&=&-\lap_x^s(x \partial_x +m+2E_x)|x|^{2s} M(x)  \\
                                  &=&-\lap_x^s(x \partial_x)|x|^{2s} M(x) -\lap_x^s(m+2E_x)|x|^{2s} M(x) \\ 
                                  &=&-\lap_x^s x( |x|^{2s}\partial_x + 2s |x|^{2s-2}x ) M(x) -(m+2E_x+4s)\lap_x^s|x|^{2s} M(x) \\ 
                                  &=&-(m+2E_x+2s)\lap_x^s|x|^{2s} M(x) \\ 
                                  &=&\frac{-2(E_x+\frac{m}{2}+s)}{A_s(E_x)}M(x)\,. 
\enar
This gives the expression
\bear
P_{M,2s+1}&=&-|x|^{2s}x \pmon \frac{A_s(E_x)}{2(E_x+\frac{m}{2}+s)}\lap_x^s \partial_x \\
   &=&-\frac{A_s(E_x-2s-1)}{2(E_x+\frac{m}{2}-s-1)}|x|^{2s} x \pmon \lap_x^s \partial_x \,,
   %&=&-\frac{A_s(E_x-2s-1)}{2(E_x+\frac{m}{2}-s-1)}x^{2s+1}  \pmon \partial_x^{2s+1}  \,.
\enar
which completes the proof.
\eop
By means of the identification (\ref{identification}) we thus obtain
\bear
\pmon=(1 - \frac{1}{\Bn +1} \Fm \Fp ) \phar
\enar
The element $\pmon$ is the extremal projector for $\gosp(1|2)$ (see also Tolstoy \cite{Tol1}, Zhelobenko \cite{Zh1}). 
For the other projections we  have under the usual identification
\bear
P_{M,2s}&=& \alpha_s(\Bn) \Bm^s \pmon \Bp^s  \\
P_{M,2s+1}&=& \frac{\alpha_s(\Bn+1)}{\Bn+s+1} \Bm^s \Fm \pmon \Bp^s \Fp  \,.
\enar

%%%%%%%%%%%%%%%%%%%%%%%%%%%%%%%%%%%%%%%%%%%%%%%%%%%%%%%%%%%%%%%%%%%%%%%%%%%%%%%%%%%%%%%%%%%%%%%%%%%%%%%%%%%%%%%%%%%%%%%%%%%%%%%%%%%%%%%%%%%%%%%%%%%%%%%%%%%%%%%%%ù

\section{\bf  ONB for spherical harmonics}

Recall the Fischer decomposition
\bear
\Pol(\R^p)=\bigoplus_{s=0}^{\infty}|u|^{2s}\Har(\R^p)=\bigoplus_{s=0}^{\infty}\bigoplus_{k=0}^{\infty}|u|^{2s}\Har_k(\R^p)\,.
\enar
Identify $\Pol(\R^m)$ with $\Pol(\R^p \oplus \R^q)$. 
Let $F(u,v) \in \Pol(\R^p \oplus \R^q)$ and apply  the Fischer decomposition in $\R^p$ and $\R^q$, then:
\bear
F(u,v)=\sum_{s,k=0}^{\infty}\sum_{r,i=0}^{\infty}|u|^{2s}|v|^{2r}H_k(u)G_i(v)
\enar
where $H_k(u) \in \Har_k(\R^p)$, $G_i(v) \in \Har_i(\R^q)$.
Project this identity on the orthogonal complement of $|x|^2 \Pol(\R^m)$. Using the relation $|x|^2=|u|^2+|v|^2$ it is clear that $P_H(|u|^2)=-P_H(|v|^2)$; therefore each power of $|v|^2$ can be replaced by a power of $-|u|^2$ in $\Pol(\R^m)/|x|^2 \Pol(\R^m)$. Hence up to some (immaterial) minus signs:
\bear
P_H(F(u,v))=\sum_{s,k,i=0}^{\infty}P_H(|u|^{2s}H_k(u)G_i(v))  \,.
\enar 
In fact, we have the following Theorem.
%%%%%%%%%%%%%%%%%%%%%%%%%%%%%%%%%%%%%%%%%%%%%%%%%%%%%%%%%%%%%%%%%%%%%%%%%%%%%%%%%%%%%%%%%%%%%%%%%%%%%%%%%%%%%%%%%%%%%%%%%%%%%%%%%%%%%%%%%%%%%%%%%%%%%%%%%%%%%%%%
\begin{theorem}\label{isomorphism_harm}
 
The map
\bear
\tau_H: \Pol(\R^p) \otimes \Har(\R^q) \ra \Har(\R^m)\,, \quad  K\otimes G \mapsto P_H(K(u)G(v))
\enar
is an $\SO(p) \times \SO(q)$-invariant isomorphism.
\end{theorem}
%%%%%%%%%%%%%%%%%%%%%%%%%%%%%%%%%%%%%%%%%%%%%%%%%%%%%%%%%%%%%%%%%%%%%%%%%%%%%%%%%%%%%%%%%%%%%%%%%%%%%%%%%%%%%%%%%%%%%%%%%%%%%%%%%%%%%%%%%%%%%%%%%%%%%%%%%%%%%%%%
\begin{bewijs}
First of all $P_H(W(u,v))=0$ iff $W(u,v)=|x|^2Q(u,v)$ for some polynomial $Q \in \Pol(\R^m)$.
Let $\{K_i(u)\}$ be a basis of $\Pol(\R^p)$ and let $\{H_j(v)\}$ be a basis of $\Har(\R^q)$. 
Suppose that $P_H(\sum_{ij} \la_{ij} K_i(u)H_j(v))=0$, then  $\sum_{ij} \la_{ij} K_i(u)H_j(v)=|x|^2Q(u,v)$ for some polynomial $Q(u,v)$ which must be harmonic in $v$. Now $\sum_{ij} \la_{ij} K_i(u)H_j(v)-|u|^2Q(u,v)=|v|^2Q(u,v)$ with the left hand side  harmonic in $v$. This is only possible if $Q=0$ and thus $\sum_{ij} \la_{ij} K_i(u)H_j(v)$ $=0$ which implies $\la_{ij}=0$.  The invariance follows easily from the $\SO(m)$-invariance of $\phar$ and the fact that 
$\Pol(\R^p) \otimes \Har(\R^q)$ is an $\SO(p) \times \SO(q)$-module. 
 \eop
\end{bewijs}
We will now  determine the harmonic projection or the map $\tau_H$ in an explicit way.
%%%%%%%%%%%%%%%%%%%%%%%%%%%%%%%%%%%%%%%%%%%%%%%%%%%%%%%%%%%%%%%%%%%%%%%%%%%%%%%%%%%%%%%%%%%%%%%%%%%%%%%%%%%%%%%%%%%%%%%%%%%%%%%%%%%%%%%%%%%%%%%%%%%%%%%%%%%%%%%%%%%%%%
\begin{theorem} \label{harmonic projection}
Let $P_k(u) \in \Har_k(\R^p)$ and  $Q_i(v) \in \Har_i(\R^q)$. Then
\bear
P_H(|u|^{2s}P_k(u)Q_i(v))=\lambda(s,k,i)|x|^{2s}P_s^{k+\frac{p-2}{2},\,i+\frac{q-2}{2}}\left(\frac{|v|^2-|u|^2}{|v|^2+|u|^2}\right)P_k(u)Q_i(v)
\enar
where the constant $\lambda(s,k,i)$ is given by
\bear
\lambda(s,k,i)=(-1)^s\bin{2s+k+i+\frac{m}{2}-2}{s}^{-1}=\frac{(-1)^s s!}{(s+k+i+\frac{m-2}{2})_s}\,.
\enar
The Fischer norm is given by
\bear
|| P_H(|u|^{2s}P_k(u)Q_i(v))||^2_m
=c(s,k,i,p,q)||P_k(u)||^2_p ||Q_i(v)||^2_q
\enar
where the constant $c(s,k,i,p,q)$ is given by
\bear
c(s,k,i,p,q):=4^s s!\frac{(k+\frac{p}{2})_s (i+\frac{q}{2})_s}{(s+k+i+\frac{m-2}{2})_s}\,.
\enar
\end{theorem}
%%%%%%%%%%%%%%%%%%%%%%%%%%%%%%%%%%%%%%%%%%%%%%%%%%%%%%%%%%%%%%%%%%%%%%%%%%%%%%%%%%%%%%%%%%%%%%%%%%%%%%%%%%%%%%%%%%%%%%%%%%%%%%%%%%%%%%%%%%%%%%%%%%%%%%%%%%%%%%%%%%%%%%
\begin{bewijs}
Put $G_{s,k,i}=|u|^{2s}P_k(u)Q_i(v)$. The Laplace operator splits into $\lap_x=\lap_u+\lap_v$ and
\bear
\lap_x^j G_{s,k,i}=\lap_u^j (|u|^{2s}P_k(u))\,Q_i(v)\,.
\enar
Consider  polar coordinates $(\rho, \xi) \in \R_+ \times S^{p-1}$ in $\R^p$: $u=\rho \xi$ with $\rho=|u|, \xi \in S^{p-1}$. The Laplace operator in $\R^p$  can be written in polar coordinates as
\bear
\lap_u=\partial_{\rho}^2+\frac{p-1}{\rho}\partial_{\rho}+\frac{1}{\rho^2}\lap_{LB,\,p}, 
\enar
where $\lap_{LB,\,p}$ denotes the Laplace-Beltrami operator on the sphere $S^{p-1}$. Spherical harmonics $P_k(u)$ in $\R^p$ are eigenfunctions of  $\lap_{LB,\,p}$ of the form $\lap_{LB,\,p}P_k(u)=-k(k+p-2)P_k(u)$.
Thus
\bear
\lap_u |u|^{2s}P_k(u)&=&\left(\partial_{\rho}^2+\frac{p-1}{\rho}\partial_{\rho}+\frac{1}{\rho^2}\lap_{LB,p}\right)\rho^{2s+k}P_k(\xi) \\
                     &=&4s(s+k+\frac{p-2}{2})\rho^{2(s-1)}P_k(u)
\enar
and recursively
\bear
\lap_u^j |u|^{2s}P_k(u)
%&=&4^j s(s-1) \ldots (s-(j-1))  \\
%                       & & \;(s+k+\frac{p-2}{2}) \ldots (s-(j-1)+k+\frac{p-2}{2})\rho^{2(s-j)}P_k(u) \\
                     =4^j (-s)_j (-s-k-\frac{p-2}{2})_j \, \rho^{2(s-j)}\, P_k(u)\,.
\enar
The extremal projector takes the form
\bear
P_H(G_{s,k,i})%&=&\sum_{j=0}^{\infty}\frac{1}{4^j j!} \frac{1}{(-E_x-\frac{m}{2}+2)_j}|x|^{2j}\lap_x^j G_{s,k,i} \\   
              &=&\sum_{j=0}^{\infty}\frac{1}{4^j j!} \frac{1}{(-2s-k-i-\frac{m}{2}+2)_j}|x|^{2j}\lap_u^j |u|^{2s}P_k(u)Q_i(v) \\
              &=&|u|^{2s}\left(\sum_{j=0}^{\infty}\frac{(-s)_j(-s-k-\frac{p-2}{2})_j}{(-2s-k-i-\frac{m}{2}+2)_j}\frac{|x|^{2j}}{|u|^{2j}}\right)\,P_k(u)Q_i(v) \\  
              &=&|u|^{2s} F\left(-s, -s-k-\frac{p-2}{2}; -2s-k-i-\frac{m}{2}+2; \frac{|x|^{2}}{|u|^{2}}\right)\,P_k(u)Q_i(v)  \,. 
\enar
The classical Jacobi polynomials can be expressed in terms of the ($_2F_1$-)hypergeometric functions by means of the relation:
\bear
P_n^{\alpha, \beta}(t)=\bin{2n+\alpha+\beta}{n}\left(\frac{t-1}{2}\right)^n F \left(-n, -n-\alpha; -2n-\alpha-\beta; \frac{2}{1-t} \right)\,.
\enar
%To make the formulas a little bit neater we will sometimes use the notation $k_p:=k+\frac{p-2}{2}$, $i_q:=i+\frac{q-2}{2}$. 
Now put $n=s$, $\alpha=k+\frac{p-2}{2}$, $\beta=i+\frac{q-2}{2}$ and $\frac{2}{1-t}=\frac{|x|^{2}}{|u|^{2}}$ or $t=\frac{|v|^2-|u|^2}{|v|^2+|u|^2}$.
To make the formulas a little bit neater we will sometimes use the notation $k_p:=k+\frac{p-2}{2}$, $i_q:=i+\frac{q-2}{2}$.
Then 
\bear
P_H(|u|^{2s}P_k(u)Q_i(v))=\lambda(s,k,i)|x|^{2s}P_s^{k_p,\,i_q}\left(\frac{|v|^2-|u|^2}{|v|^2+|u|^2}\right)P_k(u)Q_i(v)\,.
\enar 
Alternatively, the Jacobi polynomial can be expanded as follows:
\bear
P_n^{\alpha, \beta}(t)=\frac{1}{2^n}\sum_{j=0}^n \bin{\alpha+n}{j} \bin{\beta+n}{n-j}(t+1)^j (t-1)^{n-j}\,.
\enar
If $t=\frac{|v|^2-|u|^2}{|v|^2+|u|^2}$, then $t+1=\frac{2|v|^2}{|v|^2+|u|^2}$ and $t-1=\frac{-2|u|^2}{|v|^2+|u|^2}$, hence
\bear
|x|^{2s}P_s^{k_p,i_q}\left(\frac{|v|^2-|u|^2}{|v|^2+|u|^2}\right)=  
\sum_{j=0}^s \bin{k_p+s}{j} \bin{i_q+s}{s-j}(-1)^{s-j}|v|^{2j} |u|^{2(s-j)} 
\enar
and
\bea  \label{expr1}
P_H(G_{s,k,i})= \lambda(s,k,i) \left(\sum_{j=0}^s \bin{k_p+s}{j} \bin{i_q+s}{s-j}(-1)^{s-j}|v|^{2j} |u|^{2(s-j)} \right)P_k(u)Q_i(v)\,.
\ena
Consider the Fischer inner product $\inpro{\;}{\;}_m$ on $\Pol(\R^m)$. Take $F_i(u) \in \Pol(\R^p)$, $G_i(v) \in \Pol(\R^q)$. Then  
\bear
\inpro{F_1(u)G_1(v)}{F_2(u)G_2(v)}_m=\inpro{F_1(u)}{F_2(u)}_p \inpro{G_1(v)}{G_2(v)}_q \,.
\enar
Let $P_k(u) \in \Har_k(\R^p)$ and  $Q_i(v) \in \Har_i(\R^q)$ have unit norm for the Fischer inner products on $\R^p$ or $\R^q$. Then 
$|v|^{2j} |u|^{2(s-j)}P_k(u)Q_i(v)$ and $|u|^{2s}P_k(u)Q_i(v)$ are orthogonal for $ j\neq 0$ and 
\bear
|| \,|u|^{2s}P_k(u)Q_i(v) ||^2_m      
&=&||\, |u|^{2s}P_k(u) ||^2_p \, ||Q_i(v)||^2_q     \\
&=&\inpro{P_k(u)}{\lap_u^s(|u|^{2s}P_k(u))}_p \, ||Q_i(v)||^2_q         \\
&=&4^s s!(k+\frac{p}{2})_s \,|| P_k(u) ||^2_p \, ||Q_i(v)||^2_q        \\
&=&4^s s!(k+\frac{p}{2})_s \,. 
\enar  
Let us now determine $||P_H(|u|^{2s}P_k(u)Q_i(v))||^2_m$. 
There is only term in (\ref{expr1}) (corresponding to $j=0$) which contributes to the inner product above. Hence
we have that
\bear
||P_H(|u|^{2s}P_k(u)Q_i(v))||^2_m    
&=&(-1)^s \lambda(s,k,i)\bin{i_q+s}{s}||\,|u|^{2s}P_k(u)||^2_p \, ||Q_i(v)||^2_q     \\
&=&\frac{\bin{i_q+s}{s}}{\bin{2s+k+i+\frac{m}{2}-2}{s}}4^s s!(k+\frac{p}{2})_s     \\
&=&4^s s!\frac{(k+\frac{p}{2})_s (i+\frac{q}{2})_s}{(s+k+i+\frac{m-2}{2})_s} \\
&:=&c(s,k,i,p,q)\,,
\enar
which completes the proof.
\eop 
\end{bewijs}
In what follows O(N)(G)B means ortho(normal)(gonal) basis with respect to the relevant Fischer inner products. 
%%%%%%%%%%%%%%%%%%%%%%%%%%%%%%%%%%%%%%%%%%%%%%%%%%%%%%%%%%%%%%%%%%%%%%%%%%%%%%%%%%%%%%%%%%%%%%%%%%%%%%%%%%%%%%%%%%%%%%%%%%%%%%%%%%%%%%%%%%%%%%%%%%%%%%%%%%%%%%%%%

Applying Theorem \ref{isomorphism_harm} we get obviously the following result.

\begin{theorem} \label{OGB harmonics}
Let $s,k,i,n \in \N$ with $2s+k+i=n$. Let 
\bear
& &\{S_{k,l}(u),l=1,\ldots, \dim \Har_k(\R^p)\} \mbox{ be an ONB of } \Har_k(\R^p) \\
& &\{Q_{i,j}(v),j=1,\ldots, \dim \Har_i(\R^q)\} \mbox{ be an ONB of } \Har_i(\R^q) \,.
\enar
Then  
\bear
P_H(|u|^{2s}S_{k,l}(u)Q_{i,j}(v)),\; l=1,\ldots, \dim \Har_k(\R^p),\; j=1,\ldots, \dim \Har_i(\R^q)  
\enar
determine an OGB of $\Har_n(\R^m)$ 
and the basis elements are eigenfunctions of 3 commuting Laplace-Beltrami $($or Casimir$)$ operators with eigenvalues
\bear
\lap_{LB,\,p}&:&-k(k+p-2)   \\
\lap_{LB,\,q}&:&-i(i+q-2)   \\
\lap_{LB,\,m}&:&-(2s+k+i)(2s+k+i+m-2)=-n(n+m-2)   
\enar
which determine the labels $(s,k,i)$ in a unique way. 
\end{theorem}

%%%%%%%%%%%%%%%%%%%%%%%%%%%%%%%%%%%%%%%%%%%%%%%%%%%%%%%%%%%%%%%%%%%%%%%%%%%%%%%%%%%%%%%%%%%%%%%%%%%%%%%%%%%%%%%%%%%%%%%%%%%%%%%%%%%%%%%%%%%%%%%%%%%%%%%%%%%%%%%%%%%
\begin{corollary}
The decomposition
\bear
\Har_n(\R^m)=\bigoplus_{2s+k+i=n}\tau_H(|u|^{2s}\Har_k(\R^p) \otimes \Har_i(\R^q))
\enar
is multiplicity free under the action of $\SO(p) \times \SO(q)$.
\end{corollary}
\begin{opmerking}
\begin{enumerate}
\item
Consider $p=1$. Take the standard ONB $\{e_1, \ldots,e_m\}$ of $\R^m$ and the corresponding chain 
\bear
\SO(m) \supset \SO(m-1) \supset \ldots \supset \SO(2) \,.
\enar
By induction we thus obtain the Gel'fand-Zetlin basis for $\Har_n(\R^m)$.
\item
Consider $p=2$. Take the Cartan basis $\mathfrak{h}=\{L_{12}, L_{34}, \ldots,  L_{2N-1\, 2N}\}$, $N=\left[\frac{m}{2}\right]$ of $\mathfrak{so}(m)$ and the corresponding chain 
\bear
\SO(m) \supset \SO(m-2) \supset \ldots \supset \SO(2) \mbox{ or } \SO(1)
\enar
By induction we thus obtain an OGB of eigenfunctions of $\mathfrak{h}$ for $\Har_n(\R^m)$.
\item
The operators $\lap_{LB,\,p}$, $\lap_{LB,\,q}$ and $\lap_{LB,\,m}$ correspond exactly to the traditional Casimir operators coming from the $\SO(p)$, $\SO(q)$ 
and $\SO(m)$-action.
\end{enumerate}
\end{opmerking}

%%%%%%%%%%%%%%%%%%%%%%%%%%%%%%%%%%%%%%%%%%%%%%%%%%%%%%%%%%%%%%%%%%%%%%%%%%%%%%%%%%%%%%%%%%%%%%%%%%%%%%%%%%%%%%%%%%%%%%%%%%%%%%%%%%%%%%%%%%%%%%%%%%%%%%%%%%%%%%%%%%%%%

\section{\bf  ONB for spherical monogenics}

%%%%%%%%%%%%%%%%%%%%%%%%%%%%%%%%%%%%%%%%%%%%%%%%%%%%%%%%%%%%%%%%%%%%%%%%%%%%%%%%%%%%%%%%%%%%%%%%%%%%%%%%%%%%%%%%%%%%%%%%%%%%%%%%%%%%%%%%%%%%%%%%%%%%%%%%%%%%%%%%%

Consider the Fischer decomposition
\bear
\Pol(\R^p, \C_p)=\bigoplus_{s=0}^{\infty}u^{s}\Mon(\R^p, \C_p)=\bigoplus_{s=0}^{\infty}\bigoplus_{k=0}^{\infty}u^{s}\Mon_k(\R^p,\C_p)\,.
\enar
Identify $\Pol(\R^m, \C_m)$ with $\Pol(\R^p \oplus \R^q, \C_m) \cong \Pol(\R^p, \C_p) \otimes \Pol(\R^q, \C_q)$. 
Let $F(u,v) \in \Pol(\R^p \oplus \R^q,\C_m)$ and apply  the monogenic Fischer decomposition in $\R^p$ and $\R^q$, then:
\bear
F(u,v)=\sum_{s,k=0}^{\infty}\sum_{r,i=0}^{\infty}u^{s}v^{r}M_k(u)N_i(v)
\enar
where $M_k(u) \in \Mon_k(\R^p, \C_p)$, $N_i(v) \in \Mon_i(\R^q,\C_q)$.
Project this identity on the orthogonal complement of $x \Pol(\R^m,\C_m)$. Using the relation $x=u+v$ it is clear that $P_M(u)=-P_M(v)$; therefore each power of $v$ can be replaced by a power of $-u$. Hence up to some (immaterial) minus signs in the summands:
\bear
P_M(F(u,v))=\sum_{s,k,i=0}^{\infty}P_M(u^{s}M_k(u)N_i(v))
\enar 
%%%%%%%%%%%%%%%%%%%%%%%%%%%%%%%%%%%%%%%%%%%%%%%%%%%%%%%%%%%%%%%%%%%%%%%%%%%%%%%%%%%%%%%%%%%%%%%%%%%%%%%%%%%%%%%%%%%%%%%%%%%%%%%%%%%%%%%%%%%%%%%%%%%%%%%%%%%%%%%%%%%
In fact, we have the following theorem.

\begin{theorem}
The map
\bear
\tau_M: \Pol(\R^p,\C_p) \otimes \Mon(\R^q,\C_q) \ra \Mon(\R^m,\C_m)\,, \quad  K\otimes G \mapsto P_M(K(u)G(v))
\enar
is a $\Spin(p) \times \Spin(q)$-invariant isomorphism.
\end{theorem}
%%%%%%%%%%%%%%%%%%%%%%%%%%%%%%%%%%%%%%%%%%%%%%%%%%%%%%%%%%%%%%%%%%%%%%%%%%%%%%%%%%%%%%%%%%%%%%%%%%%%%%%%%%%%%%%%%%%%%%%%%%%%%%%%%%%%%%%%%%%%%%%%%%%%%%%%%%%%%%%%%%%%
\begin{bewijs} 
First of all $P_M(W(u,v))=0$ iff $W(u,v)=xQ(u,v)$ for some polynomial $Q \in \Pol(\R^m,\C_m)$.
Let $\{K_i(u)\}$ be a basis of $\Pol(\R^p,\C_p)$ and let $\{M_j(v)\}$ be a basis of $\Mon(\R^q,\C_q)$. 
Suppose that $P_M(\sum_{ij} \la_{ij} K_i(u)M_j(v))=0$, then  $\sum_{ij} \la_{ij} K_i(u)M_j(v)=xQ(u,v)$ for some polynomial $Q(u,v)$ which must be harmonic in $v$.
Thus $Q(u,v)=Q_1(u,v)+vQ_2(u,v)$ where $Q_i(u,v)$ are monogenic in $v$. Now $\partial_v(\sum_{ij} \la_{ij} K_i(u)M_j(v))=
\sum_{ij} \la_{ij} K'_i(u)\partial_vM_j(v)=0$ ($'$ denoting the main involution) and $\partial_vxQ(u,v)=0$. This yields 
\bear
0=\partial_vxQ&=&\partial_v(u+v)(Q_1+vQ_2)   \\
            &=&\partial_v(uQ_1+uvQ_2+vQ_1-|v|^2Q_2)   \\
            &=&u(m+2E_v)Q_2-(m+2E_v)Q_1-2vQ_2 \\
            &=&Q_3-2vQ_2\,.
\enar
Now $Q_2$ and $Q_3:=u(m+2E_v)Q_2-(m+2E_v)Q_1$ are monogenic in $v$, hence by the direct sum property of the (monogenic) Fischer decomposition in $v$: $Q_2=Q_3=0$ and also $Q_1=0$. Thus $Q=0$ or $\sum_{ij} \la_{ij} K_i(u)M_j(v)=0$ which implies $\la_{ij}=0$.  The invariance follows from the $\Spin(m)$-invariance of $\pmon$ and the fact that $\Pol(\R^p, \C_p) \otimes \Har(\R^q, \C_q)$ is a $\Spin(p) \times \Spin(q)$-module. 
\eop
\end{bewijs}
We will now determine the monogenic projection or the map $\tau_M$ explicitly; the formulas will be expressed in terms of the harmonic projections which were computed  explicitly in Theorem \ref{harmonic projection}.
%%%%%%%%%%%%%%%%%%%%%%%%%%%%%%%%%%%%%%%%%%%%%%%%%%%%%%%%%%%%%%%%%%%%%%%%%%%%%%%%%%%%%%%%%%%%%%%%%%%%%%%%%%%%%%%%%%%%%%%%%%%%%%%%%%%%%%%%%%%%%%%%%%%%%%%%%%%%%%%
\begin{theorem} \label{theorem9}
Let $P_k(u) \in \Mon_k(\R^p,\C_p)$, $Q_i(v) \in \Mon_i(\R^q,\C_q)$. Then
\bear
P_M(|u|^{2s}P_k(u)Q_i(v))&=&\frac{1}{2(2s+k+i)+m-2} \left((2(s+k+i)+m-2)P_H(|u|^{2s}P_k(u)Q_i(v)) \right.  \\
                         & & \; \left. -2s P_H(|u|^{2s-2}uP'_k(u)vQ_i(v)) \right)\,,  \\
P_M(u|u|^{2s}P_k(u)Q_i(v))&=&\frac{1}{2(2s+1+k+i)+m-2} \left((2s+m-p+2i)P_H(|u|^{2s}uP_k(u)Q_i(v)) \right.  \\
                         & & \; \left. -(p+2s+2k) P_H(|u|^{2s}P'_k(u)vQ_i(v)) \right)    \,.
\enar
The Fischer inner products are given by
\bear
||P_M(|u|^{2s}\,P_k(u)Q_i(v))||^2_m 
%&=&4^s s!\frac{(k+\frac{p}{2})_s (i+\frac{q}{2})_s}{(s+k+i+\frac{m-2}{2})_s} \\
&=&\frac{2(s+k+i)+m-2}{2(2s+k+i)+m-2}c(s,k,i,p,q)||P_k(u)||^2_p ||Q_i(u)||^2_q\,,  \\
||P_M(u|u|^{2s}\,P_k(u)Q_i(v))||^2_m  
%&=&4^s s!\frac{(k+\frac{p}{2})_s (i+\frac{q}{2})_s}{(s+k+i+\frac{m-2}{2})_s} \\
&=&\frac{c(s+1,k,i,p,q)}{2(s+1)} ||P_k(u)||^2_p ||Q_i(u)||^2_q \,.
\enar
\end{theorem}
%%%%%%%%%%%%%%%%%%%%%%%%%%%%%%%%%%%%%%%%%%%%%%%%%%%%%%%%%%%%%%%%%%%%%%%%%%%%%%%%%%%%%%%%%%%%%%%%%%%%%%%%%%%%%%%%%%%%%%%%%%%%%%%%%%%%%%%%%%%%%%%%%%%%%%%%%%%%%%%%%
\begin{bewijs}
The monogenic projection $P_M$ is given by
\bear
P_M=P_H \frac{x\partial_x+2E_x+m-2}{2E_x+m-2}\,. 
\enar
In the $(u,v)$-coordinate system we then have
\bear
P_M=P_H \frac{(u+v)(\partial_u + \partial_v)+2(E_u+E_v)+m-2}{2(E_u+E_v)+m-2}\,. 
\enar
Let $P_k(u) \in \Mon_k(\R^p,\C_p)$, $Q_i(v) \in \Mon_i(\R^q,\C_q)$.  Then
\bear
\partial_v (|u|^{2s}P_k(u)Q_i(v))&=&|u|^{2s}P'_k(u) \partial_v Q_i(v))=0      \\
\partial_u (|u|^{2s}P_k(u)Q_i(v))&=&2su|u|^{2s-2}P_k(u) Q_i(v))      \\
\partial_v (u|u|^{2s}P_k(u)Q_i(v))&=&-u|u|^{2s}P'_k(u) \partial_v Q_i(v)=0      \\
\partial_u (u|u|^{2s}P_k(u)Q_i(v))&=&-(p+2E_u-u\partial_u)|u|^{2s}P_k(u) Q_i(v))      \\
                                  %&=&-(p+2(2s+k))|u|^{2s}P_k(u) Q_i(v)) + 2s|u|^{2s}P_k(u) Q_i(v))    \\
                                  &=&-(p+2s+2k)|u|^{2s}P_k(u) Q_i(v))\,.    
\enar
Hence
\bear
\lefteqn{\frac{x\partial_x+2E_x+m-2}{2E_x+m-2}(|u|^{2s}P_k(u)Q_i(v)) }  \\
&=&\frac{1}{2(2s+k+i)+m-2} \left((2(s+k+i)+m-2)|u|^{2s}P_k(u)Q_i(v) \right. \\
& &\; \left.-2s |u|^{2s-2}uP'_k(u)vQ_i(v) \right)   
\enar
and
\bear
\lefteqn{\frac{x\partial_x+2E_x+m-2}{2E_x+m-2}(u|u|^{2s}P_k(u)Q_i(v)) }  \\
&=&\frac{1}{2(2s+1+k+i)+m-2} \left((2s+m-p+2i)u|u|^{2s}P_k(u)Q_i(v) \right.  \\
& & \; \left. -(p+2s+2k) |u|^{2s}P'_k(u)vQ_i(v)) \right)    \,.
\enar
Applying $\phar$ to these expressions leads to the result.
\eop
\end{bewijs}

%%%%%%%%%%%%%%%%%%%%%%%%%%%%%%%%%%%%%%%%%%%%%%%%%%%%%%%%%%%%%%%%%%%%%%%%%%%%%%%%%%%%%%%%%%%%%%%%%%%%%%%%%%%%%%%%%%%%%%%%%%%%%%%%%%%%%%%%%%%%%%%%%%%%%%%%%%%%%%%%%
\begin{theorem}  \label{theorem10}
Let $s,k,i,n \in \N$ with $s+k+i=n$. Let
\bear
& &\{P_{k,l}(u),l=1,\ldots, \dim \Mon_k(\R^p,\C_p)\} \mbox{ be an ONB of } \Mon_k(\R^p,\C_p) \\
& &\{Q_{i,j}(v),j=1,\ldots, \dim \Mon_i(\R^q,\C_q)\} \mbox{ be an ONB of } \Mon_i(\R^q,\C_q) \,.
\enar
%Let $s,k,i,n \in \N$ with $s+k+i=n$. 
Then  
\bear
P_M(u^{s}P_{k,l}(u)Q_{i,j}(v)), l=1,\ldots, \dim \Mon_k(\R^p,\C_p), j=1,\ldots, \dim \Mon_i(\R^q,\C_q)  
\enar
determine an OGB of $\Mon_n(\R^m,\R_m)$.
\end{theorem}
The results of Theorem \ref{theorem9} and \ref{theorem10} were also obtained in \cite{DSS}. The construction of their basis relies on solving a certain Vekua-type system. The solution of this system of partial differential equations yields precisely the same type of Jacobi polynomials as obtained in our approach.  

%%%%%%%%%%%%%%%%%%%%%%%%%%%%%%%%%%%%%%%%%%%%%%%%%%%%%%%%%%%%%%%%%%%%%%%%%%%%%%%%%%%%%%%%%%%%%%%%%%%%%%%%%%%%%%%%%%%%%%%%%%%%%%%%%%%%%%%%%%%%%%%%%%%%%%%%%%%%%%%%%%%
\section{\bf $\G$-operators and Scasimirs}
%%%%%%%%%%%%%%%%%%%%%%%%%%%%%%%%%%%%%%%%%%%%%%%%%%%%%%%%%%%%%%%%%%%%%%%%%%%%%%%%%%%%%%%%%%%%%%%%%%%%%%%%%%%%%%%%%%%%%%%%%%%%%%%%%%%%%%%%%%%%%%%%%%%%%%%%%
In this section we define the analogues of the Laplace-Beltrami (or Casimir) operators of Theorem \ref{OGB harmonics}. 
Define the $\G$-operators in $\R^m$, $\R^p$ and $\R^q$:
\bear
\G_x=\frac{1}{2}(\left[\partial_x,x \right]+m)\,, \quad \G_u=\frac{1}{2}(\left[\partial_u,u \right]+p)\,, 
\quad \G_v=\frac{1}{2}(\left[\partial_v,v \right]+q)   \,. 
\enar
The $\G_x$-operator does not behave as nicely as the Dirac operator $\partial_x$ does with respect to the direct sum $\R^m=\R^p \oplus \R^q$ because there is a third mixed term appearing: 
\bear
\G_x=\G_u+\G_v-(u \wedge \partial_v+v \wedge \partial_u)\,.
\enar
Define the pseudoscalars $e_P=e_1 \ldots e_p$, $e_M=e_1 \ldots e_m$, then $e_P^2=(-1)^{\frac{p(p+1)}{2}}$. Put $c=(-i)^{\frac{p(p+1)}{2}}$, then $ce_P$ is (Fischer)-self adjoint and $(ce_P)^2=1$. The projection operators
\bear
P_+=\frac{1+ce_P}{2}  \quad  \mbox{and} \quad P_-=\frac{1-ce_P}{2}
\enar 
project on the $\pm$-eigenspaces of $ce_P$. The translated $\G$-operator given by 
\bear
\G_x-\frac{m-1}{2}
\enar
plays a crucial role in Clifford analysis of one vector variable. This operator anti-commutes with $x$ and $\partial_x$ and thus anti-commutes with the odd part of $U(\gosp(1|2))$ and commutes with the even part of $U(\gosp(1|2))$. The operator in $U(\gosp(1|2)$) with this properties is up to a multiple unique and is in the  literature also known as the Scasimir (operator) and denoted by $-Sc_1$ (see e.g. \cite{ABF}). 
Define now the modified $\G$-operators or Scasimirs:
\bear
S_u=ce_P(\G_u-\frac{p-1}{2})\,, \quad S_v=ce_P(\G_v-\frac{q-1}{2})\,.  
\enar
The extra factor $ce_P$ ensures that $S_u, S_v \in \End(\Mon(\R^m))$. Remark that our definition is not symmetric because for $S_v$ one would rather expect to have the extra factor $e_Q$ instead of $e_P$. However, it seems that this definition is exactly what we need for our purposes.  
%%%%%%%%%%%%%%%%%%%%%%%%%%%%%%%%%%%%%%%%%%%%%%%%%%%%%%%%%%%%%%%%%%%%%%%%%%%%%%%%%%%%%%%%%%%%%%%%%%%%%%%%%%%%%%%%%%%%%%%%%%%%%%%%%%%%%%%%%%%%%%%%%%%%%%%%%%%%%%%%%%%%
\begin{theorem} The operators $S_u, S_v$ belong to $\End(\Mon(\R^m))$ and have the following properties$:$
\begin{itemize}
\item
$S_u, S_v$ are commuting $($Fischer$)$-self-adjoint operators and their $($anti-$)$ commutation relations with respect to $\gosp(1|2)$ follow from 
\bear
\begin{array}{ll}
S_u x=(-1)^p x S_u\,,   & S_u \partial_x=(-1)^p \partial_x S_u\\
S_v x=(-1)^{p-1} x S_v\,, & S_v \partial_x=(-1)^{p-1} \partial_x S_v   
\end{array}
\enar
\item
$S_u, S_v$ commute with the total Gamma-operator $\G_x$ and the monogenic projection $P_M$\,.
\end{itemize}
\end{theorem}
%%%%%%%%%%%%%%%%%%%%%%%%%%%%%%%%%%%%%%%%%%%%%%%%%%%%%%%%%%%%%%%%%%%%%%%%%%%%%%%%%%%%%%%%%%%%%%%%%%%%%%%%%%%%%%%%%%%%%%%%%%%%%%%%%%%%%%%%%%%%%%%%%%%%%%%%%%%%%%%%%%
\begin{bewijs}
First of all: $\displaystyle{e_P^*=(-e_p)\ldots(-e_1)=(-1)^{\frac{p(p+1)}{2}}e_P}$, hence
$(c e_P)^*=i^{\frac{p(p+1)}{2}}(-1)^{\frac{p(p+1)}{2}}e_P=ce_P$ and $S_u$ is self-adjoint because it is the product of two commuting self-adjoint operators. 
Clearly $\G_u$, $\G_v$ and $e_P$ commute, thus $S_u$, $S_v$ commute. The (anti-)commutation relations between $S_u$, $S_v$ and $u,v,\partial_u,\partial_v$ follow from
\bear
& &\{ \G_u-\frac{p-1}{2},u \}=[\G_u-\frac{p-1}{2}, v]=0  \\
& &\{ \G_u-\frac{p-1}{2},\partial_u \}=[\G_u-\frac{p-1}{2}, \partial_v]=0 
\enar
and the (anti-)commutation relations of $ce_P$ with respect to $u,\partial_u, v,\partial_v$:
\bear
& &\{ce_P, u \}= \{ce_P, \partial_u \}=\left[ce_P,v \right]=\left[ce_P,\partial_v \right]=0\,, \quad (p \mbox{ even})  \\
& &\{ce_P,v \}=\{ce_P,\partial_v \}=\left[ce_P, u \right]= \left[ce_P, \partial_u \right]=0\,, \quad (p \mbox{ odd})\,.
\enar
\end{bewijs}
\eop
By means of the projection operators $P_{\pm}$ one can decompose the space of spherical monogenics in $\R^p$ as follows:
\bear
\Mon(\R^p, \C_p)&=&P_+\Mon(\R^p, \C_p) \oplus P_-\Mon(\R^p, \C_p)   \\
                &=&\Mon^+(\R^p, \C_p) \oplus \Mon^-(\R^p, \C_p)   \,,
\enar
where $\Mon^{\pm}(\R^p, \C_p)$ are the $\pm$-eigenspaces of $ce_P \in \End(\Mon(\R^p,\C_p))$. 
%%%%%%%%%%%%%%%%%%%%%%%%%%%%%%%%%%%%%%%%%%%%%%%%%%%%%%%%%%%%%%%%%%%%%%%%%%%%%%%%%%%%%%%%%%%%%%%%%%%%%%%%%%%%%%%%%%%%%%%%%%%%%%%%%%%%%%%%%%%%%%%%%%%%%%%%%%%%%%%%%%
\begin{theorem} \label{monogenic decomposition1}
The $\Spin(p) \times \Spin(q)$-invariant building blocks of  $\Mon(\R^m, \C_m)$ are eigenfunctions of the $3$ commuting Scasimir operators $S_u,S_v$, and $\G_x-\frac{m-1}{2}$; the eigenvalues are given by 
\bear
\begin{array}{|r|r|r|r|} \hline
\mbox{building block} & S_u & S_v  &\G_x-\frac{m-1}{2} \\  \hline \hline
\phantom{u} P_M(|u|^{2s}P^+_k(u)Q_i(v))& -(k+\frac{p-1}{2})  &  -(i+\frac{q-1}{2}) &-(k+i+2s+\frac{m-1}{2})  \\  \hline
\phantom{u} P_M(|u|^{2s}P^-_k(u)Q_i(v)) &  (k+\frac{p-1}{2}) &  (i+\frac{q-1}{2}) &-(k+i+2s+\frac{m-1}{2})    \\ \hline 
P_M(u|u|^{2s}P^+_k(u)Q_i(v)) &  -(-1)^{p} (k+\frac{p-1}{2}) &  (-1)^{p}(i+\frac{q-1}{2}) &-(k+i+2s+1+\frac{m-1}{2})   \\ \hline   
P_M(u|u|^{2s}P^-_k(u)Q_i(v)) &  (-1)^p (k+\frac{p-1}{2}) &  -(-1)^{p}(i+\frac{q-1}{2})&-(k+i+2s+1+\frac{m-1}{2})  \\ \hline
\end{array} 
\enar
and determine the labels $(s,k,i)$ in a unique way. Here $P_k(u) \in \Mon_k(\R^p,\C_p)$ and $Q_i(v) \in \Mon_i(\R^q,\C_q)$. As a result, all summands in 
\bear
\Mon(\R^m,\C_m)= \bigoplus_{s,k,i \in \N} \tau_M(u^s\Mon^{\pm}_k(\R^p,\C_p) \otimes \Mon_i(\R^q,\C_q))
\enar
are necessarily Fischer orthogonal.
\end{theorem}
%%%%%%%%%%%%%%%%%%%%%%%%%%%%%%%%%%%%%%%%%%%%%%%%%%%%%%%%%%%%%%%%%%%%%%%%%%%%%%%%%%%%%%%%%%%%%%%%%%%%%%%%%%%%%%%%%%%%%%%%%%%%%%%%%%%%%%%%%%%%%%%%%%%%%%%%%%%%%%%%%%%%
\begin{bewijs}
Let $P_k(u) \in \Mon_k(\R^p,\C_p)$, $Q_i(v) \in \Mon_i(\R^q,\C_q)$. Since $S_u$ commutes with $|u|^{2s}$ and anti-commutes with $u|u|^{2s}$, it is sufficient to compute  
\bear
S_u(P^{\pm}_k(u)Q_i(v))&=&\left(ce_P(\G_u-\frac{p-1}{2})P^{\pm}_k(u)\right)Q_i(v)   \\
                       &=&-(k+\frac{p-1}{2})ce_P P^{\pm}_k(u)Q_i(v)=\mp(k+\frac{p-1}{2})P^{\pm}_k(u)Q_i(v)\,.
\enar
Now $S_v$ commutes with $|u|^{2s}$ and $S_v u=(-1)^{p-1} u S_v$; moreover $\G_v-\frac{q-1}{2} $ commutes with each $\C_p$-valued polynomial $V_k(u)$. Thus it is sufficient to determine
\bear
S_v (P^{\pm}_k(u)Q_i(v))&=&ce_P P^{\pm}_k(u) (\G_v-\frac{q-1}{2}) Q_i(v)    \\
                        &=&-(i+\frac{q-1}{2})ce_P P^{\pm}_k(u)Q_i(v)=\mp(i+\frac{q-1}{2})P^{\pm}_k(u)Q_i(v) \,.
\enar
The orthogonality follows from the fact that each summand is uniquely determined by a triple of eigenvalues of the 3 commuting  Scasimirs (which are self-adjoint).
\eop
\end{bewijs}
\begin{corollary}$(\Spin(p) \times \Spin(q)$-invariant decomposition of spherical monogenics$)$\\
%The $\Spin(p) \times \Spin(q)$-invariant label $(k,i)$ can appear at most once in the decomposition 
We have that
\bea
\Mon_n(\R^m,\C_m)=\bigoplus_{s+k+i=n}\tau_M(u^{s}\Mon_k(\R^p,\C_p) \otimes \Mon_i(\R^q,\C_q))   \label{monogenic decomposition2}  \,.
\ena
\end{corollary}
\begin{opmerking}
\begin{enumerate}
\item
Let us note that the summands in (\ref{monogenic decomposition2}) are not irreducible $\Spin(p) \times \Spin(q)$-modules because the spaces of spherical monogenics under consideration are Clifford algebra-valued and the Clifford algebra $\C_m$ is not irreducible as a $\Spin(m)$-module; it decomposes into the direct sum of (irreducible) spinor representations. It is possible to consider spinor-valued polynomials instead of Clifford algebra-valued polynomials.
In this case, however, the isomorphism between the spinor space in dimension $m$ and the tensor product
of spinor spaces in dimensions $p$ resp. $q$ is more complicated and we would need to consider
cases depending on various parities separately. 

We illustrate it, for example, in the simple case where both $p$  and $q$ are odd. In this case,
the irreducible $\Spin(m)$-module $\mS_m^{\pm}$ (any of the half-spinor representations) is isomorphic (as a $\Spin(p)\times\Spin(q)$-module) to the tensor product 
$\mS_p\otimes\mS_q,$ where $\mS_p,$ resp. $\mS_q$ are the irreducible spinor representations of the
corresponding $\Spin$ groups. Then the decomposition in the Corollary above can be regarded as an explicit form of  the branching rules under the reduction
of the symmetry from $\Spin(m)$ to $\Spin(p)\times \Spin(q).$  It can be compared with the
general results presented in \cite{Howe-bran}.

Nevertheless, if one considers Clifford algebra-valued polynomials, the 
$\Spin(p) \times \Spin(q)$-invariant label $(k,i)$ can appear at most once in this decomposition. In this sense, the above decomposition can be regarded as 
being ``\,multiplicity free\,''.  
\item
Consider $p=1$.  Let $\R^m=\R e_1 \oplus \R^{m-1}$. Then $\Mon_k(\R, \C_1)$ is non-trivial if $k=0$ and $\Mon_0(\R, \C_1) \cong \C_1 \cong a+be_1,\, a,b \in \C$.
Take the standard ONB $(e_1, \ldots,e_m)$ of $\R^m$ and the corresponding chain of subgroups 
\bear
\Spin(m) \supset \Spin(m-1) \supset \ldots \supset \Spin(2) 
\enar
where $\Spin(m-j)$ is the subgroup of $\Spin(m)$ fixing the unit vectors $e_1, \ldots, e_j$. By induction we thus obtain the Gel'fand-Zetlin basis for $\Mon_n(\R^m,\C_m)$. For $m=3$, this has also been discussed in \cite{Roman}.
\item
Consider $p=2$. Let $\R^2=\mbox{span} \{e_1, e_2\} \cong e_{12}^{\perp}$. Then $\R^m=\R^2  \oplus \R^{m-2}$.
Take the Cartan basis $\mathfrak{h}=\{M_{12}, M_{34}, \ldots, M_{2M-1\,2M}\}$ of $\mathfrak{spin}(m)$ and the corresponding chain 
$$
\Spin(m)\supset\Spin(2)\times\Spin(m-2)\supset\Spin(2)\times\Spin(2)\times\Spin(m-4)\ldots
$$
By induction we thus obtain an OGB of eigenfunctions of $\mathfrak{h}$ for $\Mon_n(\R^m,\C_m)$. We will explain this in more detail. 
\end{enumerate}
\end{opmerking}

%%%%%%%%%%%%%%%%%%%%%%%%%%%%%%%%%%%%%%%%%%%%%%%%%%%%%%%%%%%%%%%%%%%%%%%%%%%%%%%%%%%%%%%%%%%%%%%%%%%%%%%%%%%%%%%%%%%%%%%%%%%%%%%%%%%%%%%%%%%%%%%%%%%%%%%%%%%%%ù%%
\section{\bf Step two branching for $\Mon_n(\R^m, \C_m)$}
Let ${\mathbb C}_2$ be the complex Clifford algebra generated by $e_1$ and $e_2$. Define the following basis for  $\C_2$ (see also \cite{VL4}):
\bear
T_{\pm}:=\frac{1}{2}(e_1 \pm i e_2), \quad I_+:=-T_+T_-=\frac{1}{2}(1+ie_{12})=P_+, \quad
I_-:=-T_-T_+=\frac{1}{2}(1-ie_{12})=P_-\,.
\enar
The elements $T_{\pm}$ are null vectors and define the so-called Witt-basis. The elements $I_{\pm}$ are idempotents and they coincide for $p=2$ with the previously defined projections: $I_{\pm}=P_{\pm}$. Consider the vector variable $u=u_1e_1+u_2e_2 \in \R^2$. An element $P_k(u) \in \Mon_k(\R^2, \C_2)$ is of the form 
$$
(u_1-e_{12}u_2)^k(a+be_{12}+ce_1+de_2),\quad a,b,c,d \in {\mathbb C}\,.
$$
This follows immediately from the Cauchy-Kovalevska extension principle for the Dirac operator on $\R^2$. Instead of the standard basis of $\C_2$ one can also consider the null-basis. They are related as follows:
\bear
a+be_{12}+ce_1+de_2=(a-ib)I_+ + (a+ib)I_- + (c-id)T_+ + (c+id)T_-\,.
\enar
Moreover $ie_{12}I_+=I_+$, $ie_{12}T_+=T_+$, $ie_{12}I_-=-I_-$ and $ie_{12}I_-=-I_-$. Introduce complex variables in the plane $(u_1, u_2) = \R^2$:
$$
z:=u_1 + i u_2\,, \quad \bar{z}:=u_1 -i u_2\,, \quad \partial_z:=\frac{1}{2}(\partial_{u_1} - i \partial_{u_2})\,, \quad
\partial_{\bar{z}}:=\frac{1}{2}(\partial_{u_1} + i \partial_{u_2})\,.
$$
We thus obtain the following basis for $\Mon_k(\R^2,\C_2)$:
\bear
(u_1-e_{12}u_2)^k\,T_+&=&(u_1+iu_2)^k\,T_+=z^k\,T_+\\
(u_1-e_{12}u_2)^k\,I_+&=&(u_1+iu_2)^k\,P_+=z^k\,I_+\\
(u_1-e_{12}u_2)^k\,T_-&=&(u_1-iu_2)^k\,T_-=\bar{z}^k\,T_-\\
(u_1-e_{12}u_2)^k\,I_-&=&(u_1-iu_2)^k\,P_-=\bar{z}^k\,I_-\,.
\enar
It is clear that our choice of Witt-basis is compatible with the standard  complex structure in the plane spanned by $e_1$ and $e_2$.  
The powers $z^k$ and $\bar{z}^k$ are the basic holomorphic and anti-holomorphic polynomials on $\R^2$. 
Expressing the summands from Theorem \ref{monogenic decomposition1} in terms of complex coordinates induced by the Witt-basis leads to the following isomorphisms: 
\bear
|u|^{2s}\Mon^{+}_k(\R^2,\C_2)& \cong & \Span \;\bar{z}^s z^{k+s}\bina{I_+}{T_+}     \\  
|u|^{2s}\Mon^{-}_k(\R^2,\C_2)& \cong & \Span \; z^s \bar{z}^{k+s}\bina{I_-}{T_-}     \\ 
u|u|^{2s}\Mon^{+}_k(\R^2,\C_2)& \cong & \Span \; \bar{z}^s z^{k+s+1}\bina{T_-}{I_-}    \\    
u|u|^{2s}\Mon^{-}_k(\R^2,\C_2)& \cong & \Span \; z^s \bar{z}^{k+s+1}\bina{T_+}{I_+} \,.
\enar 
For $p=2$, $S_u=ie_{12}(-e_{12}L_{12}-\frac{1}{2})=iM_{12}$, thus $S_u$ is up to the imaginary unit $i$ the first element of the Cartan basis of $\mathfrak{spin}(m)$ and  $S_v=ie_{12}(\G_v-\frac{m-3}{2})$. 
In complex coordinates 
$$
S_u= i M_{12}=i\left(u_1 \partial_{u_2}-u_2 \partial_{u_1}-\frac{1}{2}e_{12}\right)=E_{\bar{z}}-E_{z}-\frac{i}{2}e_{12}\,,
$$
where $E_{z}$ and $E_{\bar{z}}$ are the Euler operators which determine the degree of homogeneity in $z$ and $\bar{z}$. 
The eigenfunctions of $S_u$ on ${\Har}(\R^2,\C_2)$ are easily found to be of the form ($k \in \N $):
\bear
S_u z^{k}\bina{I_+}{T_+}&=&-(k+\frac{1}{2})\,z^{k}\bina{I_+}{T_+} \\
S_u \bar{z}^{k}\bina{I_-}{T_-}&=&\phantom{-}(k+\frac{1}{2})\,\bar{z}^{k}\bina{I_-}{T_-}\\
S_u z^{k+1}\bina{I_-}{T_-}&=&-(k+\frac{1}{2})\,z^{k+1}\bina{I_-}{T_-} \\
S_u \bar{z}^{k+1}\bina{I_+}{T_+}&=&\phantom{-}(k+\frac{1}{2})\,\bar{z}^{k+1}\bina{I_+}{T_+}\,.
\enar
Consider the direct sum  $\R^m=\R^2  \oplus \R^{m-2}$, $x=u+v$ and let $Q_i(v) \in  {\mathcal M}_i(\R^{m-2},\C_{m-2})$. 
Then $\G_v Q_i(v)=-iQ_i(v)$  and 
\bear
S_v \bina{I_+}{T_+}Q_i(v)&=&i e_{12} \bina{I_+}{T_+}(\G_v-\frac{m-3}{2}) Q_i(v)=-(i+\frac{m-3}{2})\bina{I_+}{T_+}Q_i(v)\,\\
S_v \bina{I_-}{T_-}Q_i(v)&=&i e_{12} \bina{I_-}{T_-}(\G_v-\frac{m-3}{2}) Q_i(v)=(i+\frac{m-3}{2})\bina{I_-}{T_-}Q_i(v)\,.\\
\enar
Summarizing, we have the following characterization:
\begin{theorem} $($Step two branching$)$ \\
The space $\Mon(\R^m, \C_m)$ is  the direct sum of the $\Spin(2) \times \Spin(m-2)$-invariant summands as listed below. Each of these building blocks
occurs exactly once and they are eigenspaces of the $3$ commuting Scasimir operators 
$S_u=i M_{12}$, $S_v=ie_{12}(\G_v-\frac{m-3}{2})$ and $\G_x-\frac{m-1}{2}$. The eigenvalues are given by 
\bear
\begin{array}{|c|c|c|c|} \hline
\mbox{summand}\; (p=2)& S_u & S_v  &\G_x-\frac{m-1}{2} \\  \hline \hline
P_M(\bar{z}^s z^{k+s}\bina{I_+}{T_+}Q_i(v)) &-(k+\frac{1}{2}) & -(i+\frac{m-3}{2}) &-(k+i+2s+\frac{m-1}{2})  \\  \hline
P_M(z^s \bar{z}^{k+s}\bina{I_-}{T_-}Q_i(v)) & (k+\frac{1}{2}) &  (i+\frac{m-3}{2}) &-(k+i+2s+\frac{m-1}{2})    \\ \hline 
P_M(\bar{z}^s z^{k+s+1}\bina{T_-}{I_-}Q_i(v)) & -(k+\frac{1}{2}) & (i+\frac{m-3}{2})  &-(k+i+2s+1+\frac{m-1}{2})   \\ \hline   
P_M(z^s \bar{z}^{k+s+1}\bina{T_+}{I_+}Q_i(v)) & (k+\frac{1}{2}) &-(i+\frac{m-3}{2})&-(k+i+2s+1+\frac{m-1}{2})  \\ \hline
\end{array} 
\enar
Here $s,k,i \in \N$ and $Q_i(v) \in  {\mathcal M}_i(\R^{m-2},\C_{m-2})$.
\end{theorem}

Applying this construction for $p=2$ inductively, we thus obtain an orthogonal basis of weight vectors for ${\mathcal M}_n(\R^{m},\C_{m})$. 

%%%%%%%%%%%%%%%%%%%%%%%%%%%%%%%%%%%%%%%%%%%%%%%%%%%%%%%%%%%%%%%%%%%%%%%%%%%%%%%%%%%%%%%%%%%%%%%%%%%%%%%%%%%%%%%%%%%%%%%%%%%%%%%%%%

\subsection*{Acknowledgments}

R. L\'avi\v cka and V. Sou\v cek acknowledge the financial support from the grant GA 201/08/0397.
This work is also a part of the research plan MSM 0021620839, which is financed by the Ministry of Education of the Czech Republic.

%%%%%%%%%%%%%%%%%%%%%%%%%%%%%%%%%%%%%%%%%%%%%%%%%%%%%%%%%%%%%%%%%%%%%%%%%%%%%%%%%%%%%%%%%%%%%%%%%%%%%%%%%%%%%%%%%%%%%%%%%%%%%%%%%%%%%%%%%%

\end{document}